\documentclass[11pt]{article}

\usepackage{enumerate}
\usepackage{amsmath}
\usepackage{amssymb,latexsym}
\usepackage{amsthm}

\usepackage{anysize}
\usepackage{graphicx}
\usepackage{url}
\usepackage{color}

\newtheorem{theorem}{Theorem}[section]

\newtheorem{lemma}[theorem]{Lemma}
\newtheorem{corollary}[theorem]{Corollary}

\newtheorem{proposition}[theorem]{Proposition}

\newtheorem{remark}[theorem]{Remark}

\newcommand{\Aut}{{\rm Aut}}

\newcommand{\diag}{{\rm diag}}
\newcommand{\ord}{{\rm ord}}
\newcommand{\modd}{{\rm mod\ }}
\newcommand{\Hess}{{\rm Hess}}

\newcommand{\C}{{\bf C}}
\newcommand{\QQ}{{\bf Q}}
\newcommand{\QQQ}{\bar{\bf Q}}
\newcommand{\F}{{\bf F}}
\newcommand{\Z}{{\bf Z}}

\newcommand{\Ker}[1]{\mbox{${\rm Ker\ }{#1}$}}

\newcommand{\ga}{{\gamma}}
\newcommand{\vare}{{\varepsilon}}

\newcommand{\PP}[1]{\mbox{${\bf P}^{#1}$}}

\newcommand{\A}[1]{\mbox{${\bf A}_{#1}$}}
\newcommand{\D}[1]{\mbox{${\bf D}_{#1}$}}

\newcommand{\SSS}[1]{\mbox{${\bf S}_{#1}$}}

\begin{document}

\title{The Fermat curve $x^n+y^n+z^n$: the most symmetric non-singular
algebraic plane curve\thanks{This research was  carried out with
the support of the Italian MIUR (progetto "Strutture
Geometriche, Algebriche e Combinatoria"), and of
GNSAGA.} }
\author{Fernanda Pambianco \thanks{ Department of
Mathematics and Informatics, Perugia University, Perugia,
06123, Italy email: fernanda@dmi.unipg.it; phone: +39(075)5855006;  fax: +39(075)5855024}}

\date{}
\maketitle

\begin{abstract}
A projective non-singular plane algebraic curve of degree $d\geq 4$
is called maximally symmetric if it attains the maximum order
of the automorphism groups for complex non-singular plane
algebraic curves of degree $d$. For $d\le 7$, all such curves are known. Up to projectivities, they are the Fermat curve for $d=5,7$, see \cite{kmp1,kmp2}, the Klein quartic for $d=4$, see \cite{har}, and the Wiman sextic for $d=6$, see \cite{doi}. In this paper we work on projective plane curves defined over an algebraically closed field of characteristic zero, and we extend this result to every $d\ge 8$  showing that the Fermat curve  is the unique maximally symmetric non-singular curve of degree $d$ with $d\geq 8$, up to projectivity. For $d=11,13,17,19$, this characterization of the Fermat curve has already been obtained, see \cite{kmp2}.\\

\vspace*{0.5 cm}
\end{abstract}

\section{Introduction}

The automorphism group of an algebraic curve is one of its invariants. Apart from rational and
elliptic curves, such a group is finite. The construction and classification of curves with
large automorphism groups with respect to their genera has been considered a relevant problem
in algebraic geometry. A landmark paper in this direction is  \cite{hur} in which Hurwitz
proved his bound $|\Aut(X)|\leq 84(g(X)-1)$ valid for any complex (projective, geometrically
irreducible) algebraic curve $X$ of genus $g\geq 2$. The best known example of a curve
attaining the Hurwitz bound is the Klein quartic. Other examples are also available in the
literature but the problem of determining all such curves appears to be rather difficult; see
\cite{mach1,mach2}.

In this paper we focus on (projective) non-singular plane algebraic curves $X$ defined over an algebraically closed field $k$ of characteristic zero.
It is well known that if $\deg(X)=d\geq 4$ then the automorphism group of $X$ consists of projective maps and hence it is a subgroup of $PGL(3,k)$;
see \cite{segr} and \cite[Theorem 11.29]{hkt2008}. For this case, the Hurwitz bound reads $|\Aut(X)|\leq 42d(d-3)$. It
should be noted however that the Hurwitz bound can be attained by plane non-singular algebraic curves only for small values of $d$. This gives a
motivation for the study of non-singular plane algebraic curves which are {\em{maximally symmetric}}, that is, curves which attain the maximum
order of the automorphism groups of non-singular
plane algebraic curves of a given degree $d\geq 4$.

The main result of this paper is the following theorem:
\begin{theorem} \label{th 1} Let $X$ be a projective
non-singular plane curve of degree $d$ defined over an algebraically closed field of characteristic zero. If $d\geq 8$, then $|\Aut(X)|\leq 6d^2$ and
equality holds if and only if $X$ is projectively equivalent to the Fermat curve of equation $x^d+y^d+z^d=0$.
\end{theorem}

Theorem \ref{th 1} shows that, for $d\geq 8$, the Fermat curve is, up to a
projectivity, the unique maximally symmetric non-singular plane
curve of degree $d$. Actually, this remains true for every
$d<8$ (see \cite{kmp1,kmp2}) except for two cases, namely  $d=4$ and $d=6$.
{}From previous work, the Klein quartic and the Wiman sextic \cite{wim}
are maximally symmetric for $d=4$  and $d=6$, see \cite{har}, \cite{doi}; see also \cite{kmp0}.

Since the automorphism group of a non-singular plane curve of genus $g\ge 2$ is a finite subgroup of $PGL(3,k)$, our proof of Theorem \ref{th 1} is performed according to the well-known classification of finite subgroups of $PGL(3,k)$ which deals separately with intransitive, imprimitive, primitive and simple groups as listed in Section \ref{prem}, as types T1, T2, T3, T4. For a subgroup $G$ of $\Aut(X)$  of a curve ${\bf{V}}(f)$ satisfying the hypotheses of Theorem \ref{th 1} for $d\ge 11$, we show that if $|G|\geq 6d^2$ then ${\bf{V}}(f)$ is singular whenever $G$ is either of type T1, see Lemma \ref{lem 2.3}, or of type T2 with cyclic action on the vertices of a $G$-invariant triangle, see Lemma \ref{lem 2.5}. Moreover, if $|G|\geq 6d^2$ and $G$ is of type T2 so that $G$ induces on the vertices of a $G$-invariant triangle the symmetric group of degree $3$, then ${\bf{V}}(f)$ is either singular, or projectively equivalent to the Fermat curve of degree $d$, see Lemma \ref{lem 2.7}. Obviously, $\Aut(X)$ cannot be of type T3 or T4 since $|\Aut(X)|\geq 726$  for $d\ge 11$ while the groups of those type have order at most $360$.
The cases $8\leq d \leq 10$ require some more specific considerations for the types T1 and T2 combined with the classical Hurwitz bound.
For $d=8$ a slightly better bound is proven, namely if $\Aut(X)$ has a subgroup of order $2^7$ then $X$ is projectively equivalent to the Fermat curve, see Lemma \ref{lem 3.3}. The case $d=9$ requires longer computation but we obtain a better result, namely if $\Aut(X)$ has a subgroup of order $3^4$ then $X$ is projectively equivalent either to the Fermat curve, or the curve of equation  $x^9+y^9+z^9+\lambda x^3y^3z^3=0$ with $\lambda(\lambda^3+27)\not=0$, see Lemma \ref{lem 3.5}. It should be noted that the cases $d=8,10$ might be settled  using previous work by Aluffi and Faber \cite[Section 3.6]{aluffi}, although this would require the extension of their results to curves with non ordinary flexes.

Our paper consists of five sections and is organized as follows.

In Section \ref{prem} we collect those technical results on subgroups of $PGL(3,k)$ that play a role in our proof of Theorem \ref{th 1}. In particular, Proposition \ref{pr 2.2} is the ingredient for the proof of  Lemma \ref{lem 2.3}; similarly
Propositions \ref{pr 2.4} and \ref{pr 2.6} provide the necessary information on certain subgroups of $PGL(3,k)$ for
the proofs of  Lemmas \ref{lem 2.5} and \ref{lem 2.7}, respectively. In Section \ref{crit} we investigate plane curves $C$  of degree $d\ge 3$  that are invariant by a non-trivial element $D$ of $PGL(3,k)$ fixing the vertices of the fundamental triangle, so that $D=\diag[\varepsilon^i,\varepsilon^j,1]$ with $\varepsilon \in k^*$ and $0\leq i,j \leq v-1$ with $v=\ord(\varepsilon) < \infty$. We define a set of integers $I(D)$ depending only on $d,i,j$ and we show that if $I(D)=\emptyset$ then $C$ is singular. This criterium is used to prove
Proposition \ref{lem 1.14} that plays an important role in the proofs of the above mentioned Lemmas \ref{lem 2.3}, \ref{lem 2.5} and\ref{lem 2.7}. In Section \ref{mainpr} we state and prove these lemmas. As a corollary we obtain Theorem \ref{th 1} for $d\geq 11$. The cases $8\le d \le 10$ are investigated in Section \ref{kozott}, again the above criterium and Proposition \ref{lem 1.14} are the main ingredients, but this time we also need the classical Hurwitz formula.

An outline of the proof of Theorem \ref{th 1} appeared in \cite{pam1} and with some more details in \cite{pam2}.
Surveys on automorphism groups of algebraic curves are found in \cite{accola,far}, \cite[Chapter 11]{hkt2008}.

\section{Preliminaries}
\label{prem}

We collect those results from group theory and algebraic geometry that play a role in our proof of Theorem \ref{th 1}.
Let $\C$ denote the complex number field, $\QQQ$ \ the algebraic closure of the rational number field $\QQ$, and $k$ an algebraically
closed field of characteristic zero.
From classical results, see \cite{iit,web}, every finite subgroup $H$ of $PGL
(2,\C)$ is isomorphic to one of the following groups: \A{5}, \SSS{4}, \A{4},
\D{2\nu}($\nu\geq 2$), or $\Z_{\nu}$ ($\nu\geq1$).
Since any representation of a finite group G in $GL(n,k)$ is conjugate to a representation of the
group in $GL(n,\QQQ)$, see \cite{lan}, the classification of the finite subgroups of $PGL(n,k)$ is the same as the classification of the finite subgroups of $PGL(n,\C)$. For $n=3$, we quote it from \cite{bli}, \cite{ham} as follows.

\begin{itemize}
\item[\rm T1] Intransitive groups (the action of the group on $k^3$ is reducible):
  \begin{itemize}
  \item[\rm(1)] A diagonal abelian group of rank $\leq 2$.
  \item[\rm(2)] A group having a unique fixed point on $\PP{2}(k)$.
  \end{itemize}
\item[\rm T2] Imprimitive groups (the transitive groups with a decomposition of $k^3$ such that all
elements of the groups permute the factors as a product of transpositions):
  \begin{itemize}
  \item[\rm(3)] A group having a normal diagonal abelian subgroup $N$ such that the quotient is $G/N\cong \Z_{3}$.
  \item[\rm(4)] A group having a normal diagonal abelian subgroup $N$ such that the quotient is $G/N\cong \SSS{3}$.
  \end{itemize}
\item[\rm T3] Primitive groups (The Hessian group and its subgroups):
  \begin{itemize}
  \item[\rm(5)] A semi-direct product $\Z_{3}\times \Z_{3}\bowtie \Z_{4}$ of order 36.
  \item[\rm(6)] A semi-direct product $\Z_{3}\times \Z_{3}\bowtie \QQ_{8}$ of order 72.
  \item[\rm(7)] The Hessian group $\Z_{3}\times \Z_{3}\bowtie SL(2,\F_{3})$ of order 216.
  \end{itemize}
\item[\rm T4] Simple groups:
  \begin{itemize}
  \item[\rm(8)] The icosahedral group $\A{5}$ of order 60.
  \item[\rm(9)] The Klein group $PSL(2,\F_{7})$ of order 168.
  \item[\rm(10)] The Valentiner group $\A{6}$ of order 360.
  \end{itemize}
\end{itemize}
In particular a finite subgroup $G$ of $PGL(3,k)$ with $|G|\not\in \{36,60,72,168,216,360\}$ is of type T1 or of type T2.
Moreover, a finite subgroup $G$ of $PGL(3,k)$ of type T3 or of type T4 does not leave invariant a point, line, or triangle,
while the one of type T1 leaves invariant a point or line, and the one of type T2 leaves invariant a triangle and induces
the permutation $\Z_{3}$ or $\SSS{3}$ on the vertices, see \cite{mit}.

Let $E_i=[e_1,...,e_i], i=2,3$.

\begin{lemma}\label{lem 1.2}
Let $\alpha, \vare, \eta\in k^*$ with $\ord(\vare)=e,\,
\ord(\eta)=h,\, \alpha^e\in \langle \eta \rangle,\,
A=\diag[\alpha\vare,\alpha]$ and $C=\eta E_2$, where  $E_2=[e_1, e_2]$. Then the group $K=\langle A,C\rangle$ has order
$eh$.
\end{lemma}

\begin{proof}
We may suppose $e,h\geq 2$. Assume $\alpha^e=\eta^m$ with
$m\in [1,h]$. Let $m=m'c$ and $h=h'c$, where $c=\gcd(m,h)$.
It is easy to show that $\ord(A)=eh'=eh/c$. The map
$\chi:\langle A \rangle\times  \langle C \rangle\rightarrow K$
taking $(X,Y)$ to $XY$ is a surjective group homomorphism. It
suffices to show that $|{\rm Ker\ \chi}|=h/c$. Observe that $A^iC^j$ with
$(i,j)\in [1,eh'] \times [0,h-1]$ belongs to ${\rm Ker\ \chi}$
if and only if $i=ei'$ ($i'\in [1,h']$) and $mi'+j=0$ ($\modd
h$). There are exactly $h'=h/c$ such pairs $[i,j]$. From this the assertion follows.
\end{proof}

Let $\vare $ be a primitive $\nu$-th root of unity with $\nu\geq 2$.
Since $k$ is an algebraically closed field of characteristic $0$, any projective transformation $(A)\in PGL(3,k)$ of order
$\nu$ is conjugate to a $(D)$, where
$D=\diag[\varepsilon^h,\varepsilon^i,\varepsilon^j]$ with
$\gcd(h,i,j,\nu)=1$. Note that
$(D)=(\diag[1,\varepsilon^{i-h},\varepsilon^{j-h}])$. Thus any
cyclic subgroup of order $\nu$ in $PGL(3,k)$ must be conjugate
to either $G_{0,1}=     \langle (\diag[1,1,\vare])\rangle$ or
one of $G_{i,j}= \langle (\diag[1,\vare^i,\vare^j])\rangle$
with $1\leq i<j\leq \nu-1$ satisfying $\gcd(i,j,\nu)=1$. Obviously,
$G_{0,1}$ is conjugate to $G_{1,0}$ and $G_{1,1}$.

The following lemma shows how the cyclic
subgroups of $PGL(3,k)$ of order $\nu$ can be classified.
Any cyclic subgroup of order 2 is conjugate to $\langle \diag[-1,1,1]\rangle$. So
we assume $\nu\geq 3$.

\begin{lemma}\label{lem 1.5}
Let  $i,j,i',j'\in [1,\nu-1]$,
$i<j,i'<j'$, and $\gcd(i,j,\nu)=\gcd(i',j',\nu)=1$. Then
\begin{itemize}
\item[\rm(1)] $G_{i,j}$ is conjugate to $G_{i',j'}$ if and only if
there exists an $m\in[1,\nu-1]$ satisfying $\gcd(m,\nu)=1$ and a
permutation $\sigma\in \SSS{3}$ such that
\[
\diag[\vare_{\sigma(1)},\vare_{\sigma(2)},\vare_{\sigma(3)}]\sim
\diag[1,\vare^{i'},\vare^{j'}],
\]
where $[\vare_1,\vare_2,\vare_3]=[1,\vare^{im},\vare^{jm}]$.\\
\item[\rm(2)]
$G_{i,j}$ is conjugate to some $G_{i',j'}$ with $\gcd(i',j')=1$
and $1\leq i' < j'<\nu$. If one of $\{i',j',i'-j'\}$ is prime to
$\nu$, then $G_{i,j}$ is conjugate to $G_{1,j''}$ for some
$j''\in [2, \nu-1]$. The hypothesis holds if $\nu=p^aq^b$ for
distinct primes $p$, $q$ and nonnegative integers $a$, $b$ with $a+b>0$. \\
\item[\rm(3)] Let $\nu=p^aq^b$ for distinct primes $p$, $q$ and
nonnegative integers $a$, $b$ with $a+b>0$. Then any subgroup
$H$ of $PGL(3,k)$ isomorphic to $\Z_{\nu}$ is conjugate to one of
$G_{1,j}$ $(j\in [1,\nu-1])$.
\end{itemize}
\end{lemma}

\begin{proof} For the proof of (1) and (2) we refer to \cite{kmp2}.
(3) $H$ is conjugate to either $G_{0,1}$, which is conjugate
to $G_{1,1}$, or $G_{i,j}$ $(1\leq i<j\leq \nu-1$ with
$\gcd(i,j,\nu)=1)$. The latter group is conjugate to one of $G_{1,j'}$ with $j'\in [2,\nu-1]$ by (2).
\end{proof}

\begin{lemma}\label{lem 1.7}{\rm \cite[{\rm Lemma 1.4}]{kmp2}} Let $p\geq 5$ be a prime, $\ord(\varepsilon)=p,\,
A=\diag[1,\varepsilon,1]$, and $B=\diag[1,1,\varepsilon]$. Then any subgroup $G$ of
$PGL(3,k)$
 isomorphic to $\Z_p\times \Z_p$ is conjugate to $\langle (A),(B)\rangle$.
\end{lemma}

\begin{lemma} \label{lem 1.9}Let $G$ be a  finite subgroup of $PGL(3,k)$.
\begin{itemize}
\item[\rm(1)] If $G$ is isomorphic to $\Z_2\times \Z_4$, then it is conjugate to
$\langle(A),(B)\rangle$, where $A=\diag[1,1,-1]$ and $B=\diag[1,\varepsilon,\varepsilon]$ with $\ord(\varepsilon)=4$.\\
\item[\rm(2)] If $G$ is isomorphic to $\Z_8$, then $G$ is conjugate to
one of $\langle(A_j)\rangle$ $(j\in [1,4])$, where
$A_j=\diag[1,\varepsilon,\varepsilon^j]$ with
$\ord(\varepsilon)=8$.
\end{itemize}
\end{lemma}

\begin{proof}
For the proof of (1) we refer to the proof of LEMMA 1.6\cite{kmp2}.
(2) By Lemma  \ref{lem 1.5}, $G$ is conjugate to $\langle(A_j)\rangle$ ($j\in
[1,7]$). It is easily seen that $\langle(A_5)\rangle$, $\langle(A_6)\rangle$,
$\langle(A_7)\rangle$ are conjugate, respectively, to $\langle(A_4)\rangle$,
$\langle(A_3)\rangle$, $\langle(A_2)\rangle$. For example
$(A_5)=(\diag[\varepsilon^7,1,\varepsilon^4])$\\$=(\diag[\varepsilon,1,\varepsilon^4])^7$, so
that $\langle(A_5)\rangle$ is conjugate to $\langle(A_4)\rangle$.
\end{proof}

Let $\pi_2$ be the canonical homomorphism from $GL(2,k)$ onto
$PGL(2,k)$ such that $\pi_2(B)=(B)$. We introduce a notation to
denote an element of $GL(3,k)_{[0,0,1]}=\{A\in GL(3,k) :
A[0,0,1]=[0,0,1]\}$. For $A'=[a'_{ij}]\in GL(2,k)$ and
$a'=[a'_1,a'_2]\in k^2$, $[A',a']$ stands for the matrix
$A=[a_{ij}]\in GL(3,k)$ such that $a_{ij}=a'_{ij}$ ($i,j\in
[1,2]$), $a_{3j}=a'_j$ ($j\in [1,2]$), $a_{i3}=0$ ($i\in
[1,2]$) and $a_{33}=1$. Let $PGL(3,k)_{(0,0,1)} =\{ (A)\in
PGL(3,k) : (A)\ {\rm fixes} \ (0,0,1)\}$. Then the map
$\tau:PGL(3,k)_{(0,0,1)}\rightarrow GL(3,k)_{[0,0,1]}$ sending
$(A)$ to $A/a_{33}$ is a group isomorphism, where $A=[a_{ij}]$.
Denote by $\pi$ the map sending $[A',a']\in GL(3,k)_{[0,0,1]}$
to $A'\in GL(2,k)$, and let $\psi=\pi_2\circ \pi$. Let $G_1$
and $G_2$ be subgroups of $PGL(3,k)_{(0,0,1)}$. Then they are
conjugate in $PGL(3,k)_{(0,0,1)}$, if and only if $\tau(G_1)$
and $\tau(G_2)$ are conjugate in $GL(3,k)_{[0,0,1]}$. If they
are conjugate in $PGL(3,k)_{(0,0,1)}$, then $\psi(\tau(G_1))$
and $\psi(\tau(G_2))$ are conjugate in $PGL(2,k)$. Note that
$[E_2, a']^m=[E_2,\ ma']$. Hence $[E_2,a']$ is of finite order
if and  only if $a'=[0,0]$. Therefore if $G'$ is a finite
subgroup of $GL(3,k)_{[0,0,1]}$, then the restriction
$\pi|_{G'}$ is a injective homomorphism. In fact, if $[A',a']\in
\Ker\ \pi|_{G'}$, then $[A',a']=[E_2,a']$ is of finite order so
that $a'=[0,0]$. Consequently, $G''=G'\cap {\rm Ker\ }\psi $ is
a cyclic group $\langle [\eta E_2,a']\rangle$ with $\eta\in k^*$
and $a'\in k^2$, for $\pi(G'')=\pi|_{G'}(G'')$ is isomorphic to a finite subgroup of $k^*$ which is cyclic.

\begin{proposition} \label{pr 2.2}Assume that a finite subgroup $G_0$ of $PGL(3,k)$
fixes the point $(0,0,1)$. Then  $G_0$ is conjugate to some $G$ in
$PGL(3,k)_{(0,0,1)}$ with the following
 properties. Let $G'=\tau(G)$, $G''=G'\cap {\rm Ker\ }\psi$ and $H=\psi(G')$.
Then $G''$  is a cyclic group generated by a matrix
$\diag[\eta,\eta,1]$,
where $\ord(\eta)=|G|/|H|$.
\begin{itemize}
\item[\rm(1)] If $H$ is isomorphic to \A{5}, \SSS{4}, or \A{4},
then $|G''|>d$, provided $|G|\geq 6d^2$ with $d\geq 11$.
\item[\rm(2)] If $H$ is isomorphic to the dihedral group $\D{2\nu}\ (\nu\geq 2)$ of order $2\nu$,
then there exist $\alpha,\ \beta,\ \varepsilon,\ \eta\in k^*$ with $\ord(\varepsilon)=\nu$
and $\ord(\eta)=\mu$ such that $\alpha^\nu\in \langle \eta\rangle$, $\beta\in \{1,\sqrt{\eta}\}$ $(\sqrt{1}=1)$,
$G_d'=\langle \diag[\alpha\varepsilon,\alpha,1],\ \diag[\eta,\eta,1]\rangle$ is a subgroup of $G'$ of
order $\nu\mu$, and $G'=G_d'+G_d'[\beta e_2,\beta e_1,e_3]$ is a group of order $2\nu\mu$.
If $\mu=1$, then $\nu$ is odd. In particular $(^t G')=G'$.
\item[\rm(3)]
If $H$ is isomorphic to the cyclic group $\Z_{\nu}$ of order $\nu\geq 1$, then there exist
$\alpha,\ \ \varepsilon,\ \eta \in k^*$ with $\nu=\ord(\varepsilon)$ and $\mu=\ord(\eta)$
such that $\alpha^\nu\in \langle \eta\rangle$ and
$G'$ is a group of order $\nu\mu$ generated by $\diag[\alpha\varepsilon,\alpha,1]$ and $\diag[\eta,\eta,1]$.
In particular $(^t G')=G'$.
\end{itemize}
Let a finite subgroup $\tilde{G}$ of $PGL(3,k)$ leave the line $z$ invariant.
Then $(^t\tilde{G})$ fixes the point $(0,0,1)$. Denote $\psi(\tau((^t\tilde{G})))$ by $\tilde{H}$.
\begin{itemize}
\item[\rm(4)] If $\tilde{H}$ is isomorphic to \A{5}, \SSS{4}, or \A{4}, then
$\tilde{G}$ contains a cyclic group conjugate to $\langle(\diag[\eta,\eta,1])\rangle$ of order $\ord(\eta)>d$, provided
$d\geq 11$.
\item[\rm (5)] If $\tilde{H}$ is isomorphic to $\D{2\nu}\ (\nu\geq 2)$ or $\Z_{\nu}$, then $\tilde{G}$ is
conjugate to a finite subgroup of $PGL(3,k)$ which fixes the point $(0,0,1)$.
\end{itemize}
\end{proposition}

\begin{proof} Let $G_0'=\tau(G_0)$, $G_0''=G_0'\cap {\rm Ker\ }\psi$ and
$H_0=\psi(G_0')$. Recall that $G_0''=\langle [C',a']\rangle$ where $C'=\eta E_2$ with
$\ord(\eta)=|G_0''|=|G_0'|/|H_0|$ and $a'\in k^2$. Let $\mu=|G_0''|$. Now
${T_1}^{-1}[C',a']T_1=[C',0]$ for $T_1=[E_2,x]\in GL(3,k)_{[0,0,1]}$
such that $(\eta-1)x+a'=[0,0]$.
Define $G$ conjugate to $G_0$ by
$\tau(G)={T_1}^{-1}\tau(G_0)T_1$, so that
$G''=G'\cap {\rm Ker\ }\psi={T_1}^{-1}G_0''T_1=\langle \diag[\eta,\eta,1]\rangle$ and $H=\psi(G')=H_0$.
The finite group $H=\psi(G)$ is
known to be isomorphic to \A{5}, \SSS{4}, \A{4},
\D{2\nu}($\nu\geq 2$), or $\Z_{\nu}$($\nu\geq 1$). Clearly
$|G|=|G''||H|$. Now (1) follows, for $|H|\leq \max \{|\A{5}|,
|\SSS{4}|, |\A{4}|\}=60$.

 We prove (2). Suppose $H\cong\D{2\nu}$. Let $$H_{2\nu}=\langle
(J_2),(\diag[\varepsilon,1])\rangle,\quad
J_2=\left[\begin{array}{cc}0&1\\1&0\end{array}\right]$$
 and $\ord(\varepsilon)=\nu$. Then, since any subgroup of
 $PGL(2,k)$ isomorphic to the dihedral group $\D{2\nu}$ is conjugate to $H_{2\nu}$
(\cite[p.265]{web}), we have $(T_2)^{-1}H(T_2)=H_{2\nu}$ for
some $T_2\in GL(2,k)$. Define a new $G'$, hence $G=\tau^{-1}(G')$ as well, by
replacing $G'$ by $[T_2,[0,0]]^{-1}G'[T_2,[0,0]]$ so that $G''$
remains unchanged. Since $\psi(G')=H_{2\nu}$, there exist
$A=[A',u]$ and $B=[B',w]$ in $G'$ such that
$A'=\alpha\diag[\varepsilon,1]$, $B'=\beta J_2$ with $u,\ w\in k^2$ and
$\alpha,\ \beta\in k^*$.

We will show that $u=v=[0,0]$ may be assumed.
If $\mu>1$, then $CAC^{-1}A^{-1}=[E_2,(\eta^{-1}-1)uA'^{-1}]$ and $CBC^{-1}B^{-1}=[E_2,(\eta^{-1}-1)\beta^{-1}wJ_2]$
belong to $G'$ so that $u=w=[0,0]$.
Suppose $\mu=1$, hence $G''=\{E_3\}$ and $\psi: G'\rightarrow \psi(G')$ is a group isomorphism. For $T_3=[E_2,x]$ with $x\in k^2$,
$T_{3}^{-1}AT_{3}=[A',u']$, where $u'=[x_1(1-\alpha\varepsilon)+u_1,x_2(1-\alpha)+u_2]$. Consequently if $(1-\alpha\varepsilon)(1-\alpha)\not=0$,
then there exists $T_3$ such that $u'=[0,0]$. Even if $\alpha=1$ (resp. $\alpha\varepsilon=1$), there exists $T_3=[E_2,x]$
such that $T_3^{-1}AT_3=[A',[0,0]]$, for $u_2=0$ (resp. $u_1=0$) because $\ord(A)$ is finite. Set $T_3^{-1}BT_3=[B',w']$, where
$w'\in k^2$. Then $[B',w']^2=[\beta^2E_2,w'']\in G''$, where $w''=w'\beta J_2+w'$. Thus $[B',w']^2=E_3$, namely $\beta^2E_2=E_2$
and $w'B'+w'=[0,0]$. So $\beta=\pm 1$. In addition $\psi([B',w'])\psi([A',[0,0]])\psi([B',w'])=\psi([A',[0,0])^{-1}$,
hence $[B',w'][A',[0,0]][B',w']=[A',[0,0]]^{-1}$, namely $\diag[\alpha,\alpha\varepsilon]=\diag[(\alpha\varepsilon)^{-1},\alpha^{-1}]$ (i.e. $\alpha^2\varepsilon=1$) and $w'A'B'+w'=[0,0]$. Thus $w'A'=w'$. Since either $\alpha\not=1$ or
$\alpha\varepsilon\not=1$ (recall  $\nu\geq 2$), we have $w'=[0,0]$. Suppose $\beta=-1$ and let $T_4=\diag[-1,1,1]$.
Since $T_4^{-1}[A',[0,0]]T_4=[A',[0,0]]$ and $T_4^{-1}[B',[0,0]]T_4=[J_2,[0,0]]$, we may assume $\beta=1$ if $\mu=1$.

Since $|G''|=\mu$ and $G'/G''$ is isomorphic to $\D{2\nu}$, we have $|G'|=2\nu\mu$.
Let $A=\diag[\alpha\varepsilon,\alpha,1]$, $B=[\beta e_2,\beta e_1,e_3]$, $B'=[e_2,e_1,e_3]$, $B''=[\sqrt{\eta} e_2,\sqrt{\eta}e_1,e_3]$,
 and $C=\diag[\eta,\eta,1]$. Then $G'_d=\langle A,C\rangle$ is a subgroup of $G'$ of order $\nu\mu$ by Lemma \ref{lem 1.2}. The $2\nu\mu$-element
 sets $G_d'+G_d'B$ and $G_d'+BG_d'$ are subsets of $G'$, hence they coincide. Since $A^\nu$ and $B^2$ belong to $G''$, $\alpha^\nu$ and
 $\beta^2$ belong to $\langle \eta\rangle$. Let $\beta^2=\eta^\ell$ ($\ell\in [0,\mu-1]$). Assume first that $\ell$ is even.
 Let $\beta=\eta^{\ell/2}$. Then $B=C^{\ell/2}B'$. Hence $B'\in G'$ and $G'=G_d'+G_d'B'=G_d'+B'G_d'$. So we may assume $\beta=1$.
 Let $\beta=-\eta^{\ell/2}$ and $T_4=\diag[-1,1,1]$. Then $T_4^{-1}G'T_4=G_d'+G_d'B'=G_d'+B'G_d'$. Consequently, if $\ell$ is even,
 we may assume $\beta=1$, and $G'=G_d'+G_d'B'=G_d'+B'G_d'$. Assume next that $\ell$ is odd. According as
 $\beta=\eta^{(\ell-1)/2}\sqrt{\eta}$ or $\beta=-\eta^{(\ell-1)/2}\sqrt{\eta}$, $G'$ or $T_4^{-1}G'T_4$ coincides with
 $G_d'+G_d'B''=G_d'+B''G_d'$. So if $\ell$ is odd, we may assume $\beta=\sqrt{\eta}$ and $G'=G_d'+G_d'B''=G_d'+B''G_d'$.
 Finally assume $\mu=1$. As we have seen, $\alpha^2\varepsilon=1$ and $\alpha^\nu=1$, hence $\alpha=\varepsilon^\ell$ for
 some $\ell\in [0,\nu-1]$. Therefore $2\ell+1\equiv 0 (\mod \nu),$ hence $\nu$ is odd.

Next we prove (3). Suppose $H\cong \Z_{\nu}$ ($\nu\geq 1$). We may assume that $H$ is
generated by $(S\diag[\varepsilon,1]S^{-1})$, where $S\in GL(2,k)$ and
$\ord(\varepsilon)=\nu$. Replacing $G'$ by $[S,[0,0]]^{-1}G'[S,[0,0]]$ we may assume that
$H=\psi(G')=\langle (\diag[\varepsilon,1])\rangle$ and $G''=G'\cap \Ker \psi=\langle
C\rangle$, where $C=\diag[\eta,\eta,1]$ with $\mu=\ord(\eta)=|G''|$. Clearly $|G'|=|G''||H|=\nu\mu$.
There exist $\alpha\in k^*$ and $u\in k^2$ such
that $A=[\alpha\diag[\varepsilon,1],u]=[A',u]\in G'$. Since $A^\nu\in G''$, we have
$\alpha^\nu\in \langle \eta\rangle$.
We can argue quite similarly as in the case $H\cong\D{2\nu}$ to show that we may assume $u=[0,0]$.
$G_d'=\{A^iC^j\ :\ (i,j)\in [0,\nu-1]\times [0,\mu-1]\}$
is a subset of $G'$ of $\nu\mu$ elements, hence $G'=G_d'$.

The rest of the proof is easy: (4) follows from (1), and (5) follows from (2) and (3).
\end{proof}

\begin{proposition}\label{pr 2.4} Let $G$ be a finite subgroup of $PGL(3,k)$ permuting
cyclically the vertices of a triangle, $A=\diag[\alpha\vare,\alpha,1]$,
$C=\diag[\eta,\eta,1]$ and $E=[e_3,e_1,e_2]$ with $\nu=\ord(\varepsilon)$, $\mu=\ord(\eta)$ and
$\alpha^\nu\in \langle \eta\rangle$. Then $G$ is conjugate to a group of order $3\nu\mu$
$$\langle(A),(C),(E) \rangle=G_0+G_0(E)+G_0(E^2),$$
where $G_0=\{(A^iC^j)\ :\ i\in [0,\nu-1],\ j\in [0,\mu-1] \}$ is an abelian group of order $\nu\mu$.
\end{proposition}

\begin{proof} By our assumption, there exists a  $G$-invariant triangle
$P_1P_2P_3$ such that $G$ acts on $T=\{P_1,P_2,P_3\}$ as a cyclic permutation group of order $3$. Let $\rho: G\to \SSS{3}$ be the associated homomorphism.
Note that this cyclic group is the unique Sylow 3-subgroup of $\SSS{3}$. We
may assume $P_1=(1,0,0)$, $P_2=(0,1,0)$  and $P_3=(0,0,1)$.  Let
$G_0=\{(A')\in G;\ (A')P_3=P_3\}$ and $\rho(G)= \langle
\sigma\rangle$, where $\sigma=(1\ 3\ 2)$, namely $\sigma(P_1)=P_3$,
$\sigma(P_2)=P_1$, and $\sigma(P_3)=P_2$. There exists an $(E')\in G$
satisfying $\rho((E'))=\sigma$. $E'$ has the form
$[ae_3,be_1,ce_2]$ with $abc\not=0$. $G_0=\Ker\ \rho$, for if
$(A')\in G-\Ker\ \rho$, then $\rho(A')(P_3)=P_2$ or
$\rho(A')(P_3)=P_1$. Thus $G=\langle G_0,(E')\rangle=G_0+G_0(E')+G_0(E')^2$.
Note that if $(A')\in G_0$, then $A'$ is diagonal
so that $G_0$ is commutative. There exists a nonsingular
$S_1=\diag[a_1,a_2,a_3]$ such that ${S_1}^{-1}E'S_1=\lambda
\diag[e_3,e_1,e_2]$. So we may assume $E'=E$. In
addition, $\psi(G'_0)$, where $G'_0=\tau(G_0)$, is a finite
group isomorphic to a subgroup of $k*,$  hence isomorphic to $\Z_\nu$
($\nu\geq 1$). So, without loss of generality, $\psi(G_0')=\langle
(\diag[\varepsilon,1])\rangle.$


 Denote by $\mu$ the order of $G_0''=G_0'\cap {\rm Ker\
\psi}$. Obviously, $G_0''$ is a cyclic group generated by a $C=\diag[\eta,\eta,1]$ with
$\ord(\eta)=|G_0''|=\mu$. Moreover, there exists an $A=\diag[\alpha\varepsilon, \alpha,1]\in G_0'$
such that $\psi(A)=(\diag[\varepsilon,1])$ generates $\psi(G'_0)$. Since
$A^\nu\in \Ker{\psi}$,  $\alpha^\nu\in \langle \eta\rangle$. Clearly $G_0'=\langle
G_0'',A\rangle$ with $|G_0|=\mu\nu$ by Lemma  \ref{lem 1.2}.
\end{proof}

\begin{proposition} \label{pr 2.6}
Let $G$ be a finite subgroup of $PGL(3,k)$ inducing $\SSS{3}$ on the vertices of the
triangle $ P_1P_2P_3$ where $P_1=(1,0,0)$, $P_2=(0,1,0)$ and $P_3=(0,0,1)$.
Let $G_0$ be the isotropy subgroup of $G$ at $P_3=(0,0,1)$, $G_0'=\tau(G_0)$,
$G_0''=G_0'\cap\Ker \psi$ and $H=\psi(G_0')$. Then $G_0''$ is a cyclic group of order $\mu$
generated by $C=\diag[\eta,\eta,1]$ with $\mu=\ord(\eta)$. In addition:
\begin{itemize}
\item[\rm(1)] $H$ is isomorphic to none of \A{5}, \SSS{4}, and  \A{4}.
\item[\rm(2)] If $H$ is isomorphic to the dihedral group $\D{2\nu}\ (\nu\geq 2)$,
then there exist $\alpha,\ \beta,\ \gamma,\ \varepsilon,\ \in k^*$ with $\ord(\varepsilon)=\nu$ and $\{\alpha^\nu,\ \beta\gamma\}\subset \langle \eta\rangle$
such that up to conjugacy
$$G_0=G_d+G_d([\beta e_2,\gamma e_1,e_3]),\ \ \ G=G_0+G_0(E)+G_0(E^2),$$
where
$G_d=\langle (\diag[\alpha\varepsilon,\alpha,1]),\ (\diag[\eta,\eta,1])\rangle$ is a subgroup of $G$ of
order $\nu\mu$, and $E=[ e_3, e_1,e_2]$.
\item[\rm(3)]
If $H$ is isomorphic to the cyclic group $\Z_{\nu}$ of order $\nu\geq 1$, then $\nu=2$ and there exist
$\beta,\ \gamma \in k^*$ with $\beta\gamma\in \langle \eta\rangle$
such that up to conjugacy
$$G_0=G_d+G_d([\beta e_2,\gamma e_1,e_3]),\ \ \ G=G_0+G_0(E)+G_0(E^2),$$
where
$G_d=\langle (C)\rangle$ and $E=[ e_3, e_1,e_2]$.
\end{itemize}
\end{proposition}

\begin{proof}
We keep the notation used in the proof of Proposition \ref{pr 2.4}. By our assumption
$\rho$ is surjective. Therefore, if $(A)\in G$ then
$$A=[\delta_1 e_{\sigma(1)},\delta_2
e_{\sigma(2)},\delta_3 e_{\sigma(3)}],\,\delta_i\in k^*,\
\sigma\in \SSS{3}.$$
Clearly any $C\in G_0''$ has the form
$\diag[c,c,1]$, hence the finite group $G_0''$ is a cyclic group.
The image $H=\psi(G_0')$ is isomorphic to $\A{5}$, $\SSS{4}$, $\A{4}$,
$\Z_{\nu}\ (\nu\geq 1)$ or $\D{2\nu}\ (\nu\geq 2)$. If $(B_1),(B_2)\in G_0'$, then
$\rho((B_i)^2)=id_T$ so that $B_i^2$ are diagonal, hence
$\psi(B_1)^2$ and $\psi(B_2)^2$ commute. Therefore,  $H$
is not isomorphic to either $\A{5}$, or $\SSS{4}$, or $\A{4}$.
In fact, for $\sigma=(1\ 2\ 3)$ and $\sigma'=(2\ 3\ 4)$ in
$\A{4}$ we have $\sigma^2\sigma'^2\not=\sigma'^2\sigma^2$. Now
(1) follows.

Next assume $H\cong \D{2\nu}$. Let $s=\psi(A)$ and $t=\psi(B)$
($A,B\in G_0'$) be a pair of generators of $H$ satisfying the
defining relations of $\D{2\nu}$:  $s^\nu=t^2=(ts)^2=1$. Then
$G_0'=\langle A,B,C\rangle$, being $G_0''=\langle
C\rangle$. In fact, any $X\in G_0'$ has a $Y\in \langle
A,B,C\rangle$ such that $\psi(X)=\psi(Y)$. Since any $C'\in
G_0'$ has the form $\diag[\delta_1e_{\sigma(1)},\delta_2
e_{\sigma(2)},e_3]$ with $\sigma\in \SSS{2}$ and
$\delta_1,\delta_2\in k^*$, we see that $\ord(\psi(C'))=2$  if
$\sigma=(1\ 2)$.
We may assume $A=\diag[\alpha\vare,\alpha,1]$ and $B=[\beta
e_2,\gamma e_1,e_3]$. Note that   $s$ or $t$ must be skew-diagonal.
To be exact, if $\nu=2$, there is the second possibility
that $A=[\beta e_2,\gamma e_1,e_3]$ and $B=\diag[\alpha\vare,\alpha,1]$ so that
both $st=ts$ and $t$ are skew-diagonal.
This special case will be discussed later.
Note that $\{\alpha^\nu, \ \beta\gamma\}\subset\langle \eta\rangle$, for $A^\nu,\ B^2\in
G_0''=\langle C \rangle$.  Let
$G_{0,d}'=\langle A,\ C\rangle$, which is a subgroup of $G_0'$ such that $|G_{0,d}'|=\nu
\mu$ by Lemma \ref{lem 1.2}. Since $G_0'=\langle A,\ B,\ C\rangle$ with
$|G_{0}'|=|H||G_0''|=2\nu \mu$, we can easily verify $G_0'=G_{0,d}'+G_{0,d}'B$. Let
$G_{0,d}=\tau^{-1}(G_{0,d}')$, which equals $\langle (A),(C)\rangle$.
$G$ has an element $(E')\in PGL(3,k)$ such that
the restriction $\rho((E'))$ induces the cyclic permutation $(P_3P_2P_1)$. Clearly $E'=[\delta_1
e_3,\delta_2 e_1,\delta_3 e_2]$ for some $\delta_1,\delta_2,\delta_3 \in k^*$.
There exists a nonsingular matrix $S=\diag[a_1,a_2,a_3]$ such that
$S^{-1}E'S=\delta \diag[e_3,e_1,e_2]$ ($\delta\in k^*$). So we may assume $(E)\in G$.
Since $G_{0,d}=\Ker \rho$, we see $|G|=|\SSS{3}||G_{0,d}|=6|G_{0,d}|=6|G_{0,d}'|=6\nu \mu$. The
subset $\tilde{G}=G_0+G_0(E)+G_0(E)^2$ of $G$ coincides with $G$, because
$|\tilde{G}|=3|G_0|=6|G_{0,d}'|=|G|$.
Arguing similarly, we see that the type of the  group given by the special case (i.e. $\nu=2$, $A=[\beta e_2,\gamma e_1,e_3]$
and $B=\diag[\alpha\vare,\alpha,1]$) is the same type as the group given in the general case
(i.e. $\nu=2$, $A=\diag[\alpha\vare,\alpha,1]$ and $B=[\beta e_2,\gamma e_1,e_3]$).

Assume finally that $H\cong \Z_\nu$ and that $\psi(B)$ ($B\in G_0'$) is a generator of $H$.
Since $G$ contains an element of the form $([\beta e_2,\gamma e_1,e_3])$,  $B$ cannot be
diagonal. $B$ has the form $[\beta e_2,\gamma e_1,e_3]$, hence $\ord(\psi(B))=2$ and $\nu=2$
so that $|G_0'|=2\mu$. Note that no diagonal
$B'=[\beta e_1,\gamma e_2,e_3]\in G_0'$ satisfies $\langle\psi(B')\rangle=H$, for
$\psi(B)\in H$. Therefore $G_0'=G_0''+G_0''B$. As in the case $H\cong \D{2\nu}$ we can assume
$(E)\in G$. Let $\tilde{G}=G_0+G_0(E)+ G_0(E)^2$. This is a subset of $G$, and the
right-hand side is a disjoint sum, for $G_0=(G_0'')+(G_0'')(B)$. We claim $G\subset\tilde{G}$,
hence $G=\tilde{G}$. In fact, since $\SSS{3}=\langle(1\ 2\ 3)\rangle+(1\ 2)\langle(1\ 2\
3)\rangle$ we have $\rho(\tilde{G})=\rho(G)$. Therefore if $(X)\in G$, there exists a $(Y)\in
\tilde{G}$ such that $\rho((X))=\rho((Y))$, thus $(X)(Y)^{-1}\in G_0$, hence $(X)\in
G_0(Y)\subset \tilde{G}$, as claimed. Thus $|G|=3|G_0|=6\mu$.
\end{proof}

Let $A=[a_{ij}]\in GL(n,k)$ with $A^{-1}=[\alpha_{ij}]$, and let $k[x]=k[x_1,\dots,x_n]$. For any $f\in k[x]$,
let $T_A$ be the rational transformation of $k[x]$ defined by $T_{A}f(x)=f_A(x)$, where
$f_A(x)=f(\sum_{i=1}^n\alpha_{1i}x_i,\dots,\sum_{i=1}^n\alpha_{ni}x_i)$.
Then $T_A:k[x]\rightarrow k[x]$ is a $k$-algebra homomorphism.  Moreover, $f_{AB}=(f_B)_A$ holds for $f\in k[x]$ and
$A,B\in GL(n,k)$. Thus $T_{AB}=T_AT_B$, hence $T_A$ is a $k$-algebra isomorphism of $k[x]$.
A polynomial $f\in k[x]$ is $A$-invariant, if $f_{A^{-1}}=\lambda f$ for some $\lambda\in k^*$.
Let $f$ be an $A$-invariant non-zero polynomial such that $f_{A^{-1}}=\lambda f$, and $\nu=\ord(A)<\infty$.
Since $A^\nu=E_n$, the unit matrix, we have $\lambda^\nu=1$. The following lemma is trivial. As a result of the lemma
if $A=\diag[a_1,\dots,a_n]\in GL(n,k)$,  a homogeneous polynomial $f$ of degree $d$ satisfying $f_{A^{-1}}=\lambda f$
is the linear combination of monomials
$x_1^{i_1}\cdots x_n^{i_n}$ of degree $d$ such that $a_1^{i_1}\cdots a_n^{i_n}=\lambda$.

\begin{lemma}\label{lem 1.4}
Let $A\in GL(n,k)$ and let $f_1$, ..., $f_\ell\in k[x]$ be linearly independent homogeneous polynomials of degree $d\geq 1$ such that
$f_{jA^{-1}}=\lambda_jf_j$ $(j\in [1,\ell])$ for an $A\in GL(n,k)$. Then a
linear combination $f=c_1f_1+...+c_\ell f_\ell\not =0$ satisfies
$f_{A^{-1}}=\lambda f$ for some
$\lambda\in k$ if and only if $\Lambda=\{\lambda_i\ :\ c_i\not= 0\}$ consists of a single point $\lambda$.
\end{lemma}

The following fact is well-known.

\begin{lemma}\label{lem 1.8} Let $A=[a_{ij}]\in GL(n,k)$ and $f\in k[x_1,x_2,\dots,x_n]$ be homogeneous.
If $f$ is $A$-invariant, then the Hessian of $f$, namely $\Hess(f)=\det[f_{x_ix_j}] $, is also $A$-invariant.
\end{lemma}
\begin{lemma} \label{lem 1.11}
Let $\Hess(g)$ be the Hessian of a homogeneous polynomial $g=g(x,y)$ of degree $d\geq 1$. Then $\Hess(g)$ is the
zero polynomial if and only if there exist $[a,b]\in k^2$ such that $g\sim (ax+by)^d$.
\end{lemma}

\begin{proof}
It is clear that $\Hess((ax+by)^d)$ is identically zero. Let
$h=\Hess(g)$, and $A\in GL(2,k)$. Then $\Hess(g_{A^{-1}})=({\rm
det}\ A)^2h_{A^{-1}}$. Assuming that $d>1$ and that $g$ has
linearly independent linear factors $ax+by$ and $cx+dy$, we
will show that $h$ is not the zero polynomial. Considering
$g_A$ in  place of $g$ for some $A\in GL(2,k)$, we may assume
that $g$ is of the form $c_\ell x^\ell
y^m+c_{\ell-1}x^{\ell-1}y^{m+1}+...+c_0y^{\ell+m}$, where
$\ell, m\geq 1$ and $c_\ell\not=0$. Now, $h={c_\ell}^2 \ell
m(1-\ell-m)x^{2\ell-2}y^{2m-2}+d_{2\ell-3}
x^{2\ell-3}y^{2m-1}+...+d_0y^{2\ell+2m-4}$, which is not the
zero polynomial.
\end{proof}

Let $P\in \PP{2}(k)$, and let $f=f(x,y,z)$ and $g=g(x,y,z)$ be homogeneous polynomials. Then
$I(P,f\cap g)$ denotes the intersection number of ${\bf{V}}(f)$ and ${\bf{V}}(g)$ at $P$; see
\cite{ful}. If $f$ is irreducible, and $P$ is a non-singular point of ${\bf{V}}(f)$, then the
local ring ${\cal O}_P(f)$ is a discrete valuation ring, whose order function will be denoted
by $\ord_P^f$. For the definition of $\ord_P^f(g)$ see \cite[p.104]{ful}. We quote a theorem
\cite[p.116]{ful} and give a formula in Lemma \ref{lem 1.10}.

\begin{lemma}\label{lem 1.10}
Let $\Hess(f)$ be the Hessian of an irreducible homogeneous polynomial $f=f(x,y,z)$.
\begin{itemize}
\item[\rm(1)] $P$ lies both on ${\bf{V}}(f)$ and  ${\bf{ V}}(\Hess(f)),$ if and only if
$P$ is a flex or a multiple point of ${\bf{V}}(f)$.
\item[\rm(2)] $I(P, f \cap \Hess(f))$ is equal to
$1$ if and only if $P$ is an ordinary flex.  If $P$ is a simple point of ${\bf{V}}(f)$, and
$\ell$ is the tangent to ${\bf{V}}(f)$ at $P$, then
$$I(P, f \cap \Hess(f))={\rm ord}_P^f(\Hess(f))={\rm ord}_P^f(\ell)-2.$$
\end{itemize}
\end{lemma}

\begin{proof} It suffices to prove the equalities in (2). We may assume $P=(0,0,1)$ and $\ell={\bf V}(y)$.
Let $h=\Hess(f)(x,y,z)$, $f'=f(x,y,1)$, $h'=h(x,y,1)$, $g'=(f'_y)^2f'_{xx}+(f'_x)^2f'_{yy}-2f'_xf'_yf'_{xy}$,
and $Q=[0,0]\in k^2$. Then, by definition $I(P,f\cap h)=I(Q,f'\cap h')$. The right-hand side is equal to $I(Q,f'\cap g')$
\cite[p.116]{ful}. In addition $I(P,f\cap h)={\rm ord}_P^f(h)$ and $I(Q,f'\cap g')={\rm ord}_Q^{f'}(g')$ \cite[p.81]{ful}.
Since ${\rm ord}_P^f(\ell)={\rm ord}_Q^{f'}(y)$ by definition, it remains to show
${\rm ord}_Q^{f'}(g')={\rm ord}_Q^{f'}(y)-2 $. Let $f'=y+(c_{20}x^2+c_{21}xy+c_{22}y^2)+\cdots+(c_{d0}x^d+\cdots+c_{dd}y^d)$ with
$c_{20}=\cdots =c_{e-1,0}=0$ and $c_{e0}\not=0$ ($e\in [2,d]$). Then we can easily verify that ${\rm ord}_Q^{f'}(y)=e$, and
\[
{\rm ord}_Q^{f'}({f'_y}^2f'_{xx})=e-2,\ \ {\rm ord}_Q^{f'}({f'_x}^2f'_{yy})\geq 2e-1,\ \ {\rm ord}_Q^{f'}(f'_xf'_yf'_{xy})\geq e-1.
\]
Thus ${\rm ord}_Q^{f'}(g')=e-2$.
\end{proof}

The automorphism group of the Fermat curve is known, see \cite{leo} and \cite{shi}. In this paper we only need the following result.
\begin{lemma}\label{lem 1.13} Let $f(x,y,z)=x^d+y^d+z^d $, where $d\geq 4$ is an integer.
Then $|\Aut({\bf{V}}(f)|=6d^2$.
\end{lemma}

\section{Criterium for the singularity of a plane curve left invariant by a non-trivial linear transformation}
\label{crit}
In this section we use the following notation. For $\varepsilon \in k^*$ and $i,j\in [0,\nu-1]$,
let $\ord(\varepsilon)=\nu\geq 2$ and $D=\diag[\varepsilon^i,\varepsilon^j,1]$. Then $\ord(D)=\ord((D))$, and $\ord(D)=\nu$ with $\gcd(i,j,\nu)=1$.
Let ${\bf{V}}(f)$ be a plane algebraic curve of degree $d\geq 3$ such that $(D)$ is an automorphism of ${\bf{V}}(f)$.
Then $f=f(x,y,z)$ is a homogeneous polynomial such that
$f_{D^{-1}}=\lambda f$ with $\lambda=\varepsilon^r$ for some integer $r$. For a family $\{g_\lambda\ :\ \lambda\in \Lambda\}$
of polynomials $g_\lambda(x,y,z)$, let $\{g_\lambda\ :\ \lambda\in \Lambda\}_{A^{-1}}$ denote the family $\{g_{\lambda,A^{-1}}\ :\ \lambda\in \Lambda\}$ for
$A\in GL(3,k)$. In accordance with
\begin{eqnarray*}
\{x^d,\ x^{d-1}y,\ x^{d-1}z\}_{D^{-1}}&=&\{\varepsilon^{di}x^d,\ \varepsilon^{(d-1)i+j}x^{d-1}y,\ \varepsilon^{(d-1)i}x^{d-1}z \},\\
\{y^d,\ y^{d-1}x,\ y^{d-1}z\}_{D^{-1}}&=&\{\varepsilon^{dj}y^d,\ \varepsilon^{(d-1)j+i}y^{d-1}x,\ \varepsilon^{(d-1)j}y^{d-1}z\},\\
\{z^d,\ z^{d-1}x,\ z^{d-1}y\}_{D^{-1}}&=&\{z^d,\ \varepsilon^iz^{d-1}x,\ \varepsilon^jz^{d-1}y\},
\end{eqnarray*}
we define the following subsets $I_x(D)$, $I_y(D)$ and $I_z(D)$ of $\Z/\nu\Z$:
\[
 I_x(D)=\{di,(d-1)i+j,(d-1)i\}, I_y(D)=\{dj,(d-1)j+i,(d-1)j\}, I_z(D)=\{0,j,i\}.
\]
Their intersection
$$I(D)=I_x(D)\cap I_y(D)\cap I_z(D)$$
will be a useful tool in investigating the singularities of $D$-invariant plane curves.

\noindent Let $D'=\diag[1,\varepsilon^i,\varepsilon^j]$ and $E=[e_3,e_1,e_2]$.
Then $D'=E^{-1}DE$,
$$\{x^d,x^{d-1}y,x^{d-1}z\}_{E^{-1}}=\{y^d,y^{d-1}z,y^{d-1}x\},$$  and $$\{x^d,x^{d-1}y,x^{d-1}z\}_{E^{-2}}=\{z^d,z^{d-1}x,z^{d-1}y\}.$$
Therefore
\[
 \{x^d,x^{d-1}y,x^{d-1}z\}_{D'^{-1}}=\{x^d,\ \varepsilon^ix^{d-1}y,\ \varepsilon^j x^{d-1}z\}.
\]
So $I_x(D')=I_z(D)$ by definition.
Similarly, by definition $I_y(D')=I_x(D)$, $I_z(D')=I_y(D)$ and $I(D')=I(D)$.

\begin{lemma}\label{lem 1.6} If $f_{D^{-1}}=\varepsilon^r f$ and $r\not\in I(D)$,
then ${\bf{V}}(f)$ has a singular point. In particular if $I(D)=\emptyset$, then any
$(D)$-invariant plane curve of degree $d$ has a singular point.
\end{lemma}

\begin{proof}
By our assumption $r\in I_x(D)^c\cup I_y(D)^c\cup I_z(D)^c$.
Note that for any monomial $m=x^{d_1}y^{d_2}z^{d_3}$ we have
$m_{D^{-1}}=\varepsilon^{h}m$
 with $h=id_1+jd_2$.
If $r\in I_x(D)^c$, namely $r$ differs from $di$, $(d-1)i+j$, $(d-1)i$ ($\modd \nu$), then, by
Lemma  \ref{lem 1.4}, the coefficients of $x^d,\,x^{d-1}y$ and $x^{d-1}z$ of $f$ are all equal
to zero. Thus ${\bf{V}}(f)$ is singular at $(1,0,0)$. Similarly, ${\bf{V}}(f)$ is singular at
$(0,1,0)$ (resp. $(0,0,1)$) according as  $r\in I_y(D)^c$ or $r\in I_z(D)^c$. \end{proof}

\begin{proposition} \label{lem 1.14}
Let $A=\diag[\alpha\varepsilon,\alpha,1]$, $C=\diag[1,1,\eta]$ with $\nu=\ord(\varepsilon)$, $\mu=\ord(\eta)$,
${\alpha}^\nu \in \langle\eta\rangle$ and $f(x,y,z)$ a homogeneous polynomial of degree $d$.
\begin{itemize}
  \item[\rm(1)] If $\mu>d$ and $d\geq 2$, then any $(C)$-invariant ${\bf{V}}(f)$ is singular.
  \item[\rm(2)] If $\mu\in [2,d]$ and $\nu>d\geq 3$, then any $\{{(A),(C)}\}$-invariant ${\bf{V}}(f)$ is singular.
  \item[\rm(3)] If $\mu=1$ and $\nu\geq d^2\geq 9$, then any $(A)$-invariant ${\bf{V}}(f)$ is singular.
\end{itemize}
\end{proposition}

\begin{proof}
(1) Since $f_{C^{-1}}=\lambda f$ ($\lambda=\eta^i$ for some $i\in [0,\mu-1]$), $f$ has the form
$f_{d-i}(x,y)z^i$. If $i=0$, then $(0,0,1)$ is a singular point of ${\bf V}(f)$.
(2) Let $f_{A^{-1}}=\lambda f$ and $f_{C^{-1}}=\eta^i f$ ($i\in [0,\mu-1]$). If $i>0$, then
$z$  divides $f$. Suppose $i=0$ so that $f=\sum_{j=0}^{[d/\mu]} f_{d-\mu j}(x,y)z^{\mu j}$. We may assume
$f_{d}\not=0$. There is at most one $\ell_j\in [0,d-\mu j]$ such that $\alpha^{d-\mu j}\varepsilon^{\ell_j}=\lambda$, for
$d-\mu j<\nu$. Now the $A$-invariant $f$ has the form,
$$f=\sum_{j=0}^{[d/\mu]}c_jx^{\ell_j}y^{d-\mu j-\ell_j}z^{\mu j},$$
where $c_j=0$ ($j\geq 0$) unless there exists an $\ell_j\in [0,d-\mu j]$ such that $\alpha^{d-\mu j}\varepsilon^{\ell_j}=\lambda$.
We may assume $c_0\neq 0.$
Thus ${\bf{V}}(f)$ is singular at $(1,0,0)$ or $(0,1,0)$ according as $\ell_0\in [0,1]$ or $\ell_0\in [2,d]$.
(3) Note that $\alpha=\varepsilon^j$ ($j\in [0,\nu-1]$) and $A=D=\diag[\varepsilon^i,\varepsilon^j,1]$ with $i=j+1$. We may
assume $j\in [1,\nu-2]$ by (1). We denote $I_x(D)$, $I_y(D)$, $I_z(D)$ and $I(D)$ by $I_x$, $I_y$, $I_z$ and $I$, respectively.
Note also that $di\not\equiv dj,\,(d-1)j+i $ $(\modd\nu)$
and $(d-1)i+j\not\equiv dj,\, (d-1)j+i$ ($\modd \nu$), for $\gcd(i-j,\nu)=1$ and $r\not\equiv 0$
($\modd\nu$) if $r\in [1,d^2-1]$.
Moreover, if $pi+qj\equiv 0$ ($\modd\nu$), and $si+tj\equiv 0$
($\modd\nu$), then $(pt-qs)i\equiv 0$ ($\modd\nu$) and $(pt-qs)j\equiv 0$ ($\modd \nu$) so that
$(pt-qs)(i-j)\equiv 0$ ($\modd\nu$), hence $pt-qs\equiv 0$ ($\modd\nu$), provided $p,q,s,t$ are integers.
It suffices to show $I=\emptyset$. In the rest of the proof we omit the notation ($\modd\nu$) for simplicity.
Suppose $0\in  I_x\cap I_y$ (namely, $0\in I$). If $0\equiv di$, then we should have $0\equiv (d-1)j$.
As $(p,q,s,t)=(d,0,0,d-1)$, we get $0\equiv d(d-1)(j-i)$ so that
$d^2-d\equiv 0$, a contradiction. If $0\equiv (d-1)i+j$, then we should have $0\equiv (d-1)j$
 so that $(d-1)^2(j-i)\equiv 0$, a contradiction. If $0\equiv (d-1)i$, then either $0\equiv dj$ or $0\equiv (d-1)j+i$,
for $(d-1)i\not\equiv  (d-1)j$. Both cases are impossible, since the former one yields
$d(d-1)(j-i)\equiv 0$ while the latter one implies $(d-1)^2(j-i)\equiv 0$. Therefore
$0\not\in I$. Suppose $i\in I_x\cap I_y$.
If $i\equiv di$, then $i\equiv (d-1)j$. If $i\equiv (d-1)i+j$, then
$i\equiv (d-1)j$. If $i\equiv (d-1)i$, then either $i\equiv dj$ or $i\equiv (d-1)j+i$. The three cases imply, respectively,
$(d-1)^2\equiv 0$, $d^2-3d+3\equiv 0$, and either $d(d-2)\equiv 0$ or $(d-2)(d-1)$, a
contradiction. Similarly, $j\not\in Ix\cap I_y$. Therefore $I$ is empty.
\end{proof}

\begin{lemma}\label{lem 1.12}
In the pencil $f(x,y,z)=x^9+y^9+z^9+\lambda x^3y^3z^3$ $(\lambda \in k)$, the curve ${\bf{V}}(f)$ is
singular if and only if either $\lambda^3=-27$ $(or \;\; \lambda=\infty).$ For $\lambda(\lambda^3+27)\not=0$,
$$|\Aut({\bf{V}}(f))|=2\cdot 3^4<6\cdot 9^2.$$
\end{lemma}

\begin{proof} It is well known that ${\bf{V}}(f)$ is singular at $(a,b,c)$ if and only if $f_x$, $f_y$
and $f_z$ vanish there. A straightforward computation shows the first claim. We verify that  $\Aut({\bf{V}}(f))$ has a subgroup $N$ of
order $27$ that fixes each vertex of the fundamental triangle. In fact, $A=\diag[\alpha,\beta,1]$ defines such an automorphism if and only if $\alpha^9=\beta^9=(\alpha\beta)^3=1$, and it is easily seen that these equations have exactly $27$ common solutions. Let $G$ be the subgroup of $\Aut({\bf{V}}(f))$ that leaves the fundamental triangle invariant. Obviously, $N$ is a normal subgroup of $G$ such that $G/N$ is a subgroup of the symmetric group on the the vertices of the fundamental triangle. Therefore, $|G|\leq 2 \cdot 3^4$. Actually, this upper bound is attained since $G$ is invariant under any permutation of the indeterminates $x,y,z$. It remains to prove that $\Aut({\bf{V}}(f))=G$. If this were not true then $\Aut({\bf{V}}(f))$ would not leave invariant any point, line or triangle. From the remark
before Lemma \ref{lem 1.2}, $\Aut({\bf{V}}(f))$ would be of types T3 or T4. But this is impossible by Lagrange's theorem since no group of those types has
order a multiple of $2\cdot 3^4$.
\end{proof}

\section{The main theorem}
\label{mainpr}

In this section we prove the following result.

\begin{theorem} \label{th 2.1} Assume $d\geq 11$. Then the most symmetric non-singular
plane algebraic curve  of degree $d$ is projectively equivalent to the Fermat curve of equation $x^d+y^d+z^d=0$.
\end{theorem}
The proof is organized in a series of lemmas.
\begin{lemma} \label{lem 2.3} Let $G$ be a finite subgroup of $PGL(3,k)$ which fixes
 the point $(0,0,1)$ or leaves the line $z$ invariant. If $|G|\geq 6d^2$ and $d\geq 11$, then any $G$-invariant  plane curve  ${\bf{V}}(f)$
of degree $d$ is singular.
\end{lemma}

\begin{proof}
Assume that $G$ fixes $(0,0,1)$, and let $G'$, $G''$ and $H$ be as in Proposition \ref{pr 2.2}.
Then $|G|=|G'|=|G''||H|$.
Let $C=\diag[1,1,\eta]$ with $\ord(\eta)=\mu$ so that $G''=\langle\tau(C)\rangle$.
Suppose first that $H$ is isomorphic to either \A{5}, or \SSS{4}, or \A{4}. $G$ contains the subgroup
$\langle(C) \rangle$. Since $C^{\mu}=[e_1,e_2,e_3]$, we have $f_{C^{-1}}=\lambda f$
with $\lambda^{\mu}=1$ so that $\lambda=\eta^\ell$ for some $\ell\in [0,\mu-1]$.
By Proposition \ref{pr 2.2}(1)  $\mu>d$, hence ${\bf{V}}(f)$ is singular by Proposition \ref{lem 1.14} (1).
Suppose $H\cong \D{2\nu}$ or $H\cong \Z_{\nu}$, hence $\nu\mu\geq 3d^2$. By Proposition \ref{lem 1.14} ${\bf{V}}(f)$ is singular.
Assume finally that $G$ leaves the line $z$ invariant. Again ${\bf{V}}(f)$ is singular by Proposition \ref{pr 2.2} (4) and (5).
\end{proof}
\begin{lemma}\label{lem 2.5}
Let $G$ be a finite subgroup of $PGL(3,k)$ permuting cyclically the vertices of a triangle.
If $|G|\geq 6d^2$ and $d\geq 11$, then any $G$-invariant plane curve ${\bf{V}}(f)$ is singular.
\end{lemma}

\begin{proof}
Let $A=\diag[\alpha\varepsilon,\alpha,1]$, $C=\diag[1,1,\eta]$, and $E=[e_3,e_1,e_2]$  with $\nu=\ord(\varepsilon)$,
 $\mu=\ord(\eta)$ and $\alpha^\nu\in \langle \eta\rangle$. By Proposition \ref{pr 2.4} we may assume that
$G=G_0+G_0(E)+G_0(E^2)$ with $|G|=3\nu\mu$, where $G_0=\{(A^iC^j)\ :\ i\in [0,\nu-1],\ j\in [0,\mu-1]\}$ is a group of order $\nu\mu$.
Since $\nu\mu\geq 2d^2$, Proposition \ref{lem 1.14} immediately implies that ${\bf{V}}(f)$ is singular.
 \end{proof}
\begin{lemma} \label{lem 2.7}
Let $G$ be a finite subgroup of $PGL(3,k)$ of order at least $6d^2$ ($d\geq 11$) which induces $\SSS{3}$ on the vertices of a triangle.
Any $G$-invariant non-singular plane curve ${\bf{V}}(f)$ of degree $d$ is projectively equivalent
to the Fermat curve of equation $x^d+y^d+z^d=0$.
Moreover, $G$ has order $6d^2$, and is conjugate to
\begin{center}
$\tilde{G}=\langle(\diag[\vare,1,1]),(\diag[1,1,\vare]),([e_2,
e_1,e_3]), ([e_3, e_1, e_2]) \rangle$,
\end{center}
where $\ord(\vare)=d$. The Fermat curve of equation $x^d+y^d+z^d=0$ is the unique $\tilde{G}$-invariant non-singular plane curve of degree $d$.
\end{lemma}

\begin{proof}  We may
assume that $G$ is the group described in Proposition \ref{pr 2.6} and we use the notation
introduced there. If $H\cong \Z_2$, then $G$ contains $(C)$, where
$C=\diag[1,1,\eta]$ with $\ord(\eta)=\mu\geq d^2>d$. Hence ${\bf{V}}(f)$ is singular by Proposition \ref{lem 1.14} (1).
Assume $H\cong \D{2\nu}$. $G$ contains a subgroup $G_{0,d}=\langle (A),\ (C)\rangle$, where
$A=\diag[\alpha\vare,\alpha,1]$ and $C=\diag[1,1,\eta]$ with $\ord(\vare)=\nu$,
$\ord(\eta)=\mu$ and $\alpha^\nu\in \langle \eta\rangle$. Note that $\nu \mu=|G_{0,d}|=|G|/6\geq
d^2$. If $\mu\in [1,|G|/6]-\{d\}$, then ${\bf{V}}(f)$ is singular by Proposition \ref{lem 1.14}, a contradiction.
Therefore $\mu=d$. If $|G|>6d^2$, then $\nu>d$, hence
${\bf V}(f)$ is singular by Proposition \ref{lem 1.14} (2). So we may assume $\nu=d=\mu$ and $|G|=6d^2$.
Since $f_{C^{-1}}=\eta^i f$ with some $i\in [0,d-1]$ and ${\bf V}(f)$ is non-singular, we have $i=0$ and $f=f_d(x,y)+cz^d$ with $cf_d\not=0$.
Note that $f_{A^{-1}}=\lambda' f$ ($\lambda'\in k^*$), $({x^\ell y^{d-\ell}})_{A^{-1}}=\alpha^{d}\varepsilon^\ell x^\ell y^{d-\ell}$ ($\ell\in [0,d]$),
and $(z^d)_{A^{-1}}=z^d$. Therefore, $\lambda'=1$. If there exists only one $\ell\in [0,d]$ such that $\alpha^{d}\varepsilon^\ell=1$,
then $\ell\not\in \{0,\ d\}$, hence $\ell\in
[1,d-1]$ and ${\bf V}(f)$, where $f=bx^\ell y^{d-\ell}+cz^d$, is singular. Thus there exist $\ell,\ \ell'\in [0,d]$ such that
$\ell<\ell'$ and $\alpha^d\varepsilon^{\ell}=\alpha^d\varepsilon^{\ell'}=1$. Thus $\ell=0$, $\ell'=d$, $\alpha^d=1$ so that
$f=ax^d+by^d+cz^d$ with $a,b,c\in k^*$.
In particular ${\bf V}(f)$ is projectively equivalent to the
Fermat curve ${\bf V}(x^d+y^d+z^d)$, and $G$ is conjugate to a subgroup of $\Aut({\bf V}(x^d+y^d+z^d))$.
Since $|\Aut({\bf V}(x^d+y^d+z^d))|=6d^2$, $G$ is conjugate to $\Aut({\bf V}(x^d+y^d+z^d))$. We must show that
$\tilde{G}=\Aut({\bf V}(x^d+y^d+z^d))$. Clearly $\tilde{G}\subset \Aut({\bf V}(x^d+y^d+z^d))$, and
$G_{0,d}=\{(\diag[\varepsilon^i,\varepsilon^j,1])\ :\ i,j\in [0,d-1]\}$ is a group of order $d^2$. A matrix
$M_\sigma=[e_{\sigma(1)},e_{\sigma(2)},e_{\sigma(3)}]$ is called a permutation matrix of order $3$. As is well known,
$M_{\sigma\tau}=M_\sigma M_\tau$, and $M_{\sigma}^{-1}\diag[\varepsilon^{\nu_1},\varepsilon^{\nu_2}, \varepsilon^{\nu_3}]M_\sigma
=\diag[\varepsilon^{\nu_{\sigma(1)}},\varepsilon^{\nu_{\sigma(2)}},\varepsilon^{\nu_{\sigma(3)}}]$ for any $\sigma, \tau \in \SSS{3}$.
In particular
$(M_\sigma)^{-1}G_{0,d}(M_\sigma)=G_{0,d}$. Let $\sigma=\left(\begin{array}{ccc}1&2&3\\ 3&1&2\end{array} \right)$ and
$\tau=\left(\begin{array}{ccc}1&2&3\\ 2&1&3\end{array} \right)$. Then $E=[e_3,e_1,e_2]=M_\sigma$ and $B=[e_2,e_1,e_3]=M_\tau$.
Since $\langle \sigma\rangle \lhd \SSS{3}$, $\{M_{\gamma}\ :\ \gamma\in \SSS{3} \}=\{B^iE^j\ :\ i\in[0,1],\ j\in[0,2]\}$.
Consequently $\tilde{G}$ contains a subgroup
\[
 \sum_{\gamma\in {\rm{\bf S}}_{3}}(M_\gamma)G_{0,d}
=\sum_{\gamma\in {\rm{\bf S}}_{3}}G_{0,d}(M_\gamma)
=\sum_{\gamma\in {\rm{\bf S}}_{3}}G_{0,d}(M_\gamma)G_{0,d}
\]
of order $6d^2$. Thus $\tilde{G}=\Aut({\bf V}(x^d+y^d+z^d))$.

Let ${\bf V}(f)$ be a non-singular $\tilde{G}$-invariant plane algebraic curve of degree $d$.
The preceeding argument yields $f(x,y,z)=ax^d+by^d+cz^d$ with $abc\in k^*$. Since $f$ is $B$-invariant and $BE$-invariant, it follows that
$a=b=c$, namely ${\bf V}(f)={\bf V}(x^d+y^d+z^d)$.  \end{proof}

Theorem \ref{th 2.1} follows from the above lemmas.


\section{The case where the degree is between 8 and 10}
\label{kozott}

In this section we show that the statement of the main theorem remains valid for
$d=8,9,10.$

\begin{theorem} \label{th 3.1} For $d\in [8,10]$, the most symmetric non-singular
plane algebraic curve  of degree $d$ is projectively equivalent to the Fermat curve
of equation $x^d+y^d+z^d=0$.
\end{theorem}

 In order to prove  Theorem \ref{th 3.1} we recall that a finite group $G$ in $PGL(3,k)$ with $|G|\geq 6d^2$ ($d\geq 8$)
is of type T1 or T2, and use the following classical theorem, see (\cite{hur},\cite{kmp2}):

\begin{theorem}[Hurwitz] \label{th 3.2}
Let $X$ be a non-singular irreducible algebraic curve with genus $g=g(X).$
Assume $g'=g-1\geq 1$ and $|\Aut(X)|\geq 12g'$. Then the possible orders of $\Aut(X)$ are
\[
\hspace*{-50mm}\frac{12m}{m-6}\ g'
\]
where $m\in \{s\in \Z;\ s\geq 7\}
                 \cup \{8+\frac{4}{7},\ 16+\frac{4}{5},\ \infty \}$.
\end{theorem}

\begin{lemma} \label{lem 3.3} Let ${\bf{V}}(f)$ be a plane algebraic curve of degree $8$.
\begin{itemize}
\item[\rm(1)] Any $\Z_{5}$-invariant or $\Z_{11}$-invariant ${\bf{V}}(f)$ is singular.
\item[\rm(2)]Let $G$ be a subgroup of $PGL(3,k)$ with $|G|=2^7$. Any $G$-invariant  non-singular
${\bf{V}}(f)$ is projectively equivalent to the Fermat curve of equation $x^8+y^8+z^8=0$.
\end{itemize}
\end{lemma}

\begin{proof}
(1) Let $\ord(\varepsilon)=5$, $\ord(\eta)=11$,
$A_j=\diag[\varepsilon, \varepsilon^j,1]$ ($j\in [1,4]$), and
$B_j=\diag[\eta,\ \eta^j,1]$ ($j\in [1,10]$). Then
 any subgroup $G$ isomorphic to $\Z_5$ (resp. $\Z_{11}$) is conjugate to
$\langle (A)\rangle$ (resp. $\langle (B)\rangle$) by Lemma \ref{lem 1.5} (3). We can easily
verify that $I(A_j)=\emptyset$ and $I(B_j)=\emptyset$. To skip some computations note that
$\langle (A_j)\rangle$ ($j\in \{2,3,4\}$) are conjugate, that $\langle (B_j)\rangle$ ($j\in \{3,4,5,7,8,9\}$)
are conjugate and that $\langle (B_j)\rangle$ ($j\in \{2,6,10\}$) are conjugate. Now ${\bf{V}}(f)$ is singular by Lemma
\ref{lem 1.6}.

 (2) Since $2^7\not\in \{ 36,\ 60, \ 72,\ 168,\ 216,\ 360\}$, $G$ fixes a point, a line or a triangle, see \cite{mit}. Since $|G|$ is not a multiple of
$3$, $G$ leaves invariant a point or a line by Propositions \ref{pr 2.2}, \ref{pr 2.4} and
\ref{pr 2.6}. So we may assume that $G$ satisfies the condition (2) or (3) of Proposition
\ref{pr 2.2}, even if $G$ fixes a line. By this proposition $G$ contains a subgroup
$$G_d=\langle(\diag[\alpha\varepsilon,\alpha,1]),(\diag[\eta,\eta,1])\rangle$$ of order
$\nu\mu \in \{2^6,\ 2^7\}$. Here $\nu=\ord(\varepsilon)$, $\mu=\ord(\eta)$, and $\alpha^\nu\in
\langle\eta\rangle$. Let $A=\diag[\alpha\varepsilon,\alpha,1]$, $C=\diag[1,1,\eta]$. Observe that
$\langle (C)\rangle=\langle (\diag[\eta,\eta,1])\rangle$. If $\mu\geq 2^4$ or $\mu\in \{2^2,\ 2,\
1\}$, then ${\bf{V}}(f)$ is singular by Proposition \ref{lem 1.14} (1), (2) and (3).

Thus $\mu=2^3$ and $\nu=2^3$ or $\nu=2^4$. Since $\ord(C)=\mu$ and
$f_{C^{-1}}=\lambda f$, we have $\lambda=\eta^i$ ($i\in [0,\mu-1]$). Unless $i=0$, $z$ divides
$f$. So $i=0$, and $f=f_8(x,y)+cz^8$ ($c\in k$). Note that $f_8c \not=0$, for ${\bf{V}}(f)$ is
non-singular. The condition $f_{A^{-1}}=\xi f$ ($\xi\in k^*$) yields $\xi=1$, and
$f_{8A^{-1}}=f_8$, namely $\alpha^8f_8(\varepsilon x,y)=f_8(x,y)$. Thus $f_8$ is a monomial in $x,y$, hence, ${\bf{V}}(f)$ is
singular, provided $\nu=2^4$. Therefore $\nu=2^3$.
Since $(x^jy^{8-j})_{A^{-1}}=\alpha^8\varepsilon^j (x^jy^{8-j})$, $\alpha^8f_8(\varepsilon x,y)=f_8(x,y)$,  and $f_8$ cannot be a monomial,
we have $\alpha^8=1$ and $f_8(x,y)=ax^8+by^8$ with $ab\not=0$, namely $f=ax^8+by^8+cz^8$.
Therefore ${\bf{V}}(f)$ is projectively equivalent to the Fermat curve ${\bf V}(x^8+y^8+z^8)$.
\end{proof}

\begin{corollary} \label{cor 3.4}
Theorem \ref{th 3.1} holds for $d=8$.
\end{corollary}

\begin{proof} In this case, $g'=g-1=d(d-3)/2=4\cdot 5$. The possible values of
$|\Aut({\bf{V}}(f))|$ greater than $6d^2=2^7\cdot 3=(19+1/5)g'$ are $84g'$, $48g'$, $40g'$, $36g'$, $30g'$, $(26+2/5)g'$, $24g'$, $21g'$, and $20g'$. They are multiples of $5$ or $11$. Theorem \ref{th 3.1} for $d=8$ holds by Lemma
\ref{lem 3.3} and Lemma \ref{lem 1.13}.\end{proof}

\begin{remark} For curves with only ordinary flexes, Corollary \ref{cor 3.4} also follows from \cite[Section 3.6]{aluffi}.
\end{remark}

\begin{lemma} \label{lem 3.5} Let $G$ be a finite subgroup of $PGL(3,k)$, ${\bf{V}}(f)$ a $G$-invariant plane algebraic curve of degree $d=9$.
\begin{itemize}
\item[\rm (1)] If $G\cong\Z_5$ or $G\cong\Z_2\times \Z_4$, then ${\bf{V}}(f)$ is singular.
\item[\rm (2)] If $G\cong\Z_8$ and ${\bf{V}}(f)$ is non-singular,
then $|\Aut({\bf{V}}(f))|\leq 6\cdot 9^2=2\cdot 3^5$.
\item[\rm (3)] If $|G|=3^4$ and ${\bf{V}}(f)$ is non-singular, then ${\bf{V}}(f)$ is projectively equivalent
either to the Fermat curve so that $|\Aut({\bf{V}}(f))|=6d^2$ or to the curve of equation $x^9+y^9+z^9+\lambda x^3y^3z^3=0$ with
$\lambda(\lambda^3+27)\not=0$ so that $|\Aut({\bf{V}}(f))|=2d^2$.
 \end{itemize}
\end{lemma}

\begin{proof}
(1) Let $G\cong \Z_5$. We may assume $G=\langle (A_j)\rangle$ where
$A_j=\diag[1,\varepsilon,\varepsilon^j]$ with $\ord(\varepsilon)=5$ and $j\in [1, 4]$ by Lemma \ref{lem 1.5} (3).
It is easily seen that $I(A_j)=\emptyset$. Next let $G\cong \Z_2\times \Z_4$. By Lemma
\ref{lem 1.9} (1) we may assume $G=\langle (A),(B)\rangle$, where
$A=\diag[1,\varepsilon,\varepsilon]$ and $B=\diag[1,1,\varepsilon^2]$ with
$\ord(\varepsilon)=4$. Now, $f_{A^{-1}}=\varepsilon^i f$ for some $i\in [0,3]$, and
 $f_{B^{-1}}=(-1)^jf$ for some $j\in [0,1]$. Since $f$ is a linear
combination of monomials $x^{9-p-q}y^pz^q$ with $p+q=i$ ($\modd
4$), $f$ is divisible by $x$ unless $i=1$. Similarly, unless
$j=0$, $f$ is divisible by $z$. Even if $i=1$ and $j=0$, $f$ is
divisible by $y$.

 (2) Let $\ord(\varepsilon)=8$ and
$A_j=\diag[1,\varepsilon,\varepsilon^j]$. We may assume $G=\langle(A_j)\rangle$ for with $j\in
[1,4]$ by Lemma \ref{lem 1.9}. Since $I(A_j)=\{0,1,j\}$, ${\bf{V}}(f)$ is singular unless
$f_{A_j^{-1}}=\varepsilon^i f$ with $i\in I(A_j)$. \\
2-1) Let ${\bf{V}}(f)$ be an $(A_1)$-invariant curve. If $f_{A_1^{-1}}=f$, then $f=x^9f_0+xf_8(y,z)$, hence $x$
divides $f$. Suppose $f_{A_1^{-1}}=\varepsilon f$. Then
$$f=x^8(ay+bz)+\sum_{p=0}^9 c_py^{9-p}z^p$$
with $ay+bz\not=0$. Notice that any $T=[t_{ij}]\in GL(3,k)$ with $t_{1j}=t_{j1}=0$
($j\in[2,3]$) commutes with $A_1.$ Thus $(f_T)_{A_1^{-1}}=(f_{A_1^{-1}})_T=\varepsilon f_T$.
Obviously ${\bf{V}}(f)$ is non-singular if and only if so is ${\bf V}(f_T)$.
There exists an $S\in GL(3,k))$ whose first, second and third row is $[1,0,0]$, $[0,a,b]$ and $[0,a',b']$, respectively.
Let $T=S$.  Considering $f_T$,
we can assume $f=x^8y+g(y,z)$, where
$g=\sum_{p=0}^9 F_py^{9-p}z^p$. Since $F_9\not=0$, we may further assume $F_0=0$ (substitute
$y$ and $z-ey$ by $y$ and $z$ respectively, where $g(y,ey)=0$). Observe that $F_1\not=0$. Otherwise,
${\bf{V}}(f)$ is singular at $(0,1,0)$.   Simple computation yields $\Hess(f)=8x^6h(x,y,z)$,
where $h=7y\Hess(g)- 8x^8g_{zz}$. If $\Hess(g)=0$, then $g(y,z)=(ay+bz)^9$ by Lemma \ref{lem
1.11}, so that ${\bf{V}}(f)$ is singular. Thus $\Hess(g)\not=0$. Since the coefficient of
$x^0$ in $h$, namely $7y\Hess(g)$ is not equal to zero, $x$ is a linear factor of multiplicity
$6$ of $\Hess(f)$. Let ${\cal{L}}_6$ be the set of linear factors of $\Hess(f)$ of
multiplicity $6$. Clearly $|{\cal{L}}_6|\leq 3$, for $\deg h= 15$. Let $(B)\in
{\Aut({\bf{V}}(f))}_x$. We may assume $[b_{11},b_{12},b_{13}]=[1,0,0]$ for the first row of $B=[b_{ij}]$.
Let $Y=b_{22}y+b_{23}z$ and $Z=b_{32}y+b_{33}z$. Since
$f_{B^{-1}}\sim f$, it follows that $b_{21}=b_{31}=b_{23}=0$. In fact, $f_{B^{-1}}=\sum
b_{i_1,i_2,i_3}x^{i_1}Y^{i_2}Z^{i_3}$ with $b_{i_1,i_2,i_3}=0$ ($i_1\in [1,9]-\{8\}$), and $b_{1,8,0}=F_1b_{31}$.
Thus $b_{31}=0$. Now, the condition
$$f_{B^{-1}}=x^8(b_{21}x+b_{22}y+b_{23}z)+\sum_{p=1}^{9}F_p(b_{21}x+Y)^{9-p}Z^p\sim x^8y+g(y,z)$$ yields
$b_{21}=0$, hence $b_{23}=0$ as well. In addition $f_{B^{-1}}=b_{22}f$. Let $a=b_{22}$, $b=b_{32}$, and $c=b_{33}$.
Denote by $M$ the $(2,2)$-matrix
whose first and second row are $[a,0]$ and $[b,c]$, respectively. Now, the rows of $B$ are
$[1,0,0]$, $[0,a,0]$ and $[0,b,c]$. Since $f_{B^{-1}}=(x^8y)_{B^{-1}}+g_{B^{-1}}=af$ and
$g_{B^{-1}}=g_{M^{-1}}$,  we have $g_{M^{-1}}=ag$. In particular
$a^9g(1,b/a)=g_{M^{-1}}(1,0)=ag(1,0)=0$, for $g$ is homogeneous and $g(1,0)=0$ thanks to
$F_0=0$. By comparison of coefficients of $z^9$ and $yz^8$, we get $c^9=a$ and
$(F_8+9bF_9/a)c^8=F_8$. If $b=0$, the coefficient of $y^8z$ gives $a^8c=a$, hence $d^{64}=1$
and $a=c^9$. Thus the number of $B$ with $b=0$ is bounded by $64$. Suppose $b\not=0$. Then
$F_8\not=0$ by the equality $(F_8+9bF_9/a)c^8=F_8$, and a possible value of $c$ is a solution
to $(F_8+9F_9\ga)c^8=F_8$, where $\ga$ is a nonzero solution to $g(1,z)=0$ such that
$F_8+9F_9\ga\not=0$. The number of such $\ga$ does not exceed $8$, for $g(1,0)=0$. Each $\ga$
gives $8$ $c$'s, hence it gives $8$ 3-tuples $(c,a,b)$, where $a=c^9$ and $b=\ga a$.
Therefore the number of  $B$ with $b\not=0$ is bounded by $64$, that is,
$|\Aut({\bf{V}}(f))_x|\leq 2\cdot 8^2$. Thus
$|\Aut({\bf{V}}(f))|\leq |{\cal{L}}_6|\ |\Aut({\bf{V}}(f))_x| \leq 6\cdot 8^2<6\cdot 9^2$.\\
2-2) Let $f_{A_2^{-1}}=\varepsilon^i f$ with $i\in \{0,1,2\}=I(A_2)$. If $i=1$, $f$ is a
linear combination of monomials $x^{9-p-q}y^pz^q$ satisfying
 $p+2q=1$ ($\modd 8$), so that $y$ divides $f$, for $p$ is odd. By the involution
$B=[e_3,e_2,e_1]$ we have $B^{-1}A_2B=BA_2B=\vare^2A_2^{-1}$.
Hence, if $i=2$, then $f=(f_{A_2^{-1}})_{A_2}=\varepsilon^2
f_{A_2}$, so that
\[f_B=\varepsilon^2 (f_{A_2})_{B}=\varepsilon^2 f_{BA_2}=\varepsilon^2
f_{BA_2B^2}=\varepsilon^2(f_B)_{BA_2B}=\varepsilon^2 (f_{B})_{\varepsilon^2A_2^{-1}}
=(f_B)_{A_2^{-1}}.\]
Observe that $\Aut({\bf{V}}(f))$ and $\Aut({\bf{V}}(f_B))$ are conjugate and that ${\bf{V}}(f)$ is non-singular if and only if
 ${\bf{V}}(f_B)$ is non-singular.  So it suffices to consider the case $i=0$. Thus
$f$ is a linear combination of monomials $x^9$, $x^5z^4$, $x^4y^2z^3$, $x^3y^4z^2$, $x^2y^6z$,
$xy^8$, $xz^8$, and $y^2z^7$. If $f$ does not contain the monomial $x^9$, then ${\bf{V}}(f)$ is singular
at $(1,0,0)$. If $f$ does not contain
 either the monomial $xy^8$ or $xz^8$, then ${\bf{V}}(f)$ is singular either at $(0,1,0)$
  or $(0,0,1)$. Thus we may assume that
\begin{center}
$f=x^9+F_5x^5z^4+F_4x^4y^2z^3+F_3x^3y^4z^2+F_2x^2y^6z+x(y^8+z^8)+F_0y^2z^7$,
\end{center}
where $F_0\not=0$, for otherwise $x$ divides $f$. Consider the affine curve with equation
\begin{eqnarray*}
 -x&=&x^9+F_5x^5z^4+F_4x^4z^3+F_3x^3z^2+F_2x^2z+xz^8+F_0z^7\\
   &=&F_2zx^2+F_3z^2x^3+(F_0z^7+F_4z^3x^4)+(F_5z^4x^5+x^9).
\end{eqnarray*}
Obviously  $x$ is the tangent at $[0,0]$ which corresponds to $Q=(0,1,0)$ in the projective
plane. In addition $z$ is a uniformizing parameter at $[0,0]$ \cite[p.70]{ful}. In other
words, denoting the discrete valuation of the local ring ${\cal O}_{[0,0]}$ of the curve
${\bf{V}}(f(x,1,z))$ by $\ord_Q^f$, we have $\ord_Q^f(z)=1$. Since $\ord_Q^f(x)\geq 2,$ see \cite[p.71]{ful},
and $$\ord_Q^f(a_1+\cdots +a_m)\geq \min\{\ord_Q^f(a_1),\dots,\ord_Q^f(a_m)\},$$ see
\cite[p.48]{ful}, we have $\ord_Q^f(x)=\ord_Q^f(-x)\geq 5$. Therefore, evaluating again, we
obtain $\ord_Q^f(x)=\ord_Q^f(F_0z^7)=7$, for $F_0\not=0$ and the equality $$\ord_Q^f(a_1+\cdots
+a_m)= \min\{\ord_Q^f(a_1),\dots,\ord_Q^f(a_m)\}$$ holds provided that there is an $i\in [1,m]$
such that $\ord_Q^f(a_i)<\ord_Q^f(a_j)$ for any $j\not=i,$ see \cite[p.48]{ful}. By Lemma \ref{lem
1.10} the intersection number $r$ at $Q$ is equal to
 $$I(Q,f \cap \Hess(f)) =\ord_Q^f(x)-2= 5.$$
  Let ${\cal P}_r=\{P\ ;\ I(P, f \cap \Hess(f))=r\}$, whose size is not greater than
$3d(d-2)/r$ by B\'ezout's theorem. Since $\Aut({\bf{V}}(f))$ acts on ${\cal P}_r$, see
\cite[p.22]{lan},
 ${\cal P}_5\supset \ \Aut({\bf{V}}(f))Q,$
so that
 $|\Aut({\bf{V}}(f))|\leq |\Aut({\bf{V}}(f))_Q|3d(d-2)/5$. Let $(B)\in \Aut({\bf{V}}(f))_Q$.
Since $(B)$ fixes the tangent  $x$ as well as $Q=(0,1,0)$, $B$ can be assumed to have the first row
$[1,0,0]$, the second row $[b_{21},b_{22},b_{23}]$, and the third row $[b_{31},0,b_{33}]$. Clearly
$f_{B^{-1}}=\lambda f$ for some $\lambda\in k^*$. Writing $f_{B^{-1}}=\sum_{i=0}^9h_{9-i}(x,y)z^i$, we get
$b_{23}=0$, for $h_0=0$. Now $h_2(x,y)z^7\sim y^2z^7$ implies $b_{21}=b_{31}=0$, for $b_{22}b_{33}\not=0$. So $B=\diag[1,b,c]$ and
$f_{B^{-1}}=f$, hence $b^8=c^8=1$ and $b^2c^7=1$, namely $b^2=c$. Thus $|\Aut({\bf{V}}(f))_Q|\leq 8$ so that
$|\Aut({\bf{V}}(f))|<6d^2$. \\
2-3) Suppose $f_{A_3^{-1}}=\varepsilon^i f$ for some $i\in \{0,1,3\}=I(A_3)$. If $i=0$, then
$f$ is a linear combination of monomials $x^{9-p-q}y^pz^q$ such that $p+3q\equiv 0$ ($\modd 8$),
hence $p+q\equiv 0$ ($\modd 2$), so that $x$ divides $f$. The involution $B=[e_1,e_3,e_2]$ satisfies
 $B^{-1}A_3 B=A_3^{3}$. If $i=3$, then $\varepsilon
f=f_{A_3^{-3}}=(f_{BA_3^{-1}B})=(f_B)_{BA_3^{-1}}$ so that $\varepsilon f_B=(f_B)_{A_3^{-1}}$.
Therefore it suffices to consider the case $i=1$. Thus $f$ is a linear combination of
monomials
$$x^8y,x^6z^3,x^4y^3z^2, x^2y^6z, x^2y^2z^5, y^9, y^5z^4,yz^8.$$
If the coefficient of $y^9$ is equal to zero, then ${\bf{V}}(f)$ is singular at $(0,1,0)$. If
either $yx^8$ or $yz^8$ is missing in $f$, then ${\bf{V}}(f)$ is singular at $(1,0,0)$ or
$(0,0,1)$. Now, a non-singular ${\bf{V}}(f)$ such that $f_{A_3^{-1}}=\varepsilon f$ is
projectively equivalent to an $f'=f_{D^{-1}}$ with $D=\diag[d_1,d_2,d_3]$ of the form
\begin{center}
$f'=y^9+F_6x^2y^6z+F_5y^5z^4+F_3x^4y^3z^2+F_2x^2y^2z^5+y(x^8+z^8)+F_0x^6z^3$,
\end{center}
where $F_0\not=0$ (Observe that $A_3D=DA_3$ so that $f'_{A_3^{-1}}=\varepsilon f'$). We denote
$f'$ by $f$ again. $R=(0,0,1)$ is a flex of ${\bf{V}}(f)$, with tangent  $y$.
 Denote the discrete valuation of
the local ring ${\cal O}_{[0,0])}$ of the affine curve ${\bf{V}}(f(x,y,1))$ by $\ord_R^f$. In this
case, with similar evaluation as in the case 2-2), $\ord_R^f(y)=6$, so that
 $I(R, f \cap \Hess(f))=4$. Let $(B)\in \Aut({\bf{V}}(f))_R$. Since $(B)$ fixes the line $y$ and the point
$R$, the first, the second and the third row of $B$ may be $[b_{11},b_{12},0]$, $[0,1,0]$ and
$[b_{31},b_{32},b_{33}]$. Writing $f_{B^{-1}}=\sum_{i=0}^9h_{9-i}(x,y)z^i$, we get $b_{31}=b_{32}=b_{12}=0$,
for the condition $f_{B^{-1}}\sim f$ yields $h_0(x,y)=0$ and $h_6(x,y)z^3\sim x^6z^3$. Now $B=\diag[a,1,c]$ and
$f_{B^{-1}}=f$, hence $a^8=c^8=a^6c^3=1$ so that $c=a^6$, $a^8=1$. Thus
$|\Aut({\bf{V}}(f))|\leq
|\Aut({\bf{V}}(f))_R||{\cal P}_4|\leq 6d(d-2)$, where
${\cal P}_4=\{P\ :\ I(P,f \cap \Hess(f))=4 \}$.\\
2-4) Finally suppose $f_{A_4^{-1}}=\varepsilon^i f$ for some
$i\in \{0,1,4\}=I(A_4)$. If $i=1$, $f$ is divisible by $y$, for
$f$ is a linear combination of monomials $x^{9-p-q}y^pz^q$
satisfying $p+4q\equiv 1$ ($\modd 8$), hence $p\equiv 1$ ($\modd 2$). The
involution $B=[e_3,e_2,e_1]$ satisfies
$B^{-1}A_4B=\vare^4A_4^{5}$. Hence, if $i=4$, then
$-f=f_{A_4^{-5}}=f_{-BA_4^{-1}B} =(f_B)_{-BA_4^{-1}}$, hence
$f_B=(f_B)_{A_4^{-1}}$. So it suffices to consider the case
$i=0$ alone. Thus $f$ is a linear combination of monomials
$x^9$, $x^7z^2$, $x^5z^4$, $x^4y^4z$, $x^3z^6$, $x^2y^4z^3$,
$xy^8$, $xz^8$ and $y^4z^5$. Without loss of generality
\begin{center}
$f=x^9+F_7x^7z^2+F_5x^5z^4+F_4x^4y^4z+F_3x^3z^6+F_2x^2y^4z^3+x(y^8+z^8)+F_0y^4z^5$,
\end{center}
where $F_0\not=0$. $Q=(0,1,0)$ is a flex of ${\bf{V}}(f)$, with tangent $x$. Since $F_0\not=0$, we have
$\ord_Q^f(x)=5$ so that $I(Q, f \cap \Hess(f))=3$. Similarly $R=(0,0,1)$ is a flex of ${\bf{V}}(f)$, with
tangent $x$, and $\ord_R^f(x)=4$ so that $I(R, f \cap \Hess(f))=2$. A similar argument to that used in  the
preceding case 2-2) shows that any $(B)\in \Aut({\bf{V}}(f))_Q$ has the form $(\diag[1,b,c])$ with
$c=b^4$ and $b^8=1$. Conversely, if $B=\diag[1,b,c]$ with $b^8=1$ and $c=b^4$, then $(B)\in
\Aut({\bf{V}}(f))_Q$. Thus $|\Aut({\bf{V}}(f))_Q|=8$. Since $3|{\cal P}_3| +2|{\cal P}_2|\leq 3d(d-2)$ by  B\'ezout's
theorem, we have $|{\cal P}_3|<d(d-2)$. As $\Aut({\bf{V}}(f))$ acts on ${\cal P}_3$, we obtain
 $|\Aut({\bf{V}}(f))|\leq |\Aut({\bf{V}}(f))_Q||{\cal P}_3|<8d(d-2)$. As will be mentioned in the proof
of Corollary \ref{cor 3.6}, $8d(d-2)=56g'/3$ is equal to
$$\min \{|\Aut({\bf{V}}(f'))|>6d^2\ :\  {\bf{V}}(f')\ {\rm is\ a\ non-singular\ plane\ algebraic\ curve\ of\ degree}\  d \},$$
where  $g'=d(d-3)/2=3^3$.
Consequently $|\Aut({\bf{V}}(f))|\leq 6d^2=18g'$.

 (3) Since $|G|=3^4\not\in\{36,\ 60,\ 72,\ 168,\ 216,\ 360\}$, $G$ is conjugate to one of the following groups
$K_i$ ($i=1,2$) by Propositions \ref{pr 2.2}, \ref{pr 2.4} and \ref{pr 2.6}, even if $G$ fixes a line.
\begin{eqnarray*}
K_1&=&\langle (\diag[\alpha\varepsilon,\alpha,1]),\ (\diag[\eta,\eta,1])\rangle \ \  {\rm with\ } \ |K_1|=\nu\mu=3^4.\\
K_2&=&\langle (\diag[\alpha\varepsilon,\alpha,1]),\ (\diag[\eta,\eta,1]),
\ ([e_3,e_1,e_2])\rangle \ \ {\rm with\ } \ |K_2|=3\nu\mu=3^4.
\end{eqnarray*}
Observe that the case Proposition \ref{pr 2.2} (2) and the case of Proposition \ref{pr 2.6} are impossible, for $|G|$ is odd. Here
$\ord(\varepsilon)=\nu,\ \ord(\eta)=\mu$, and $\alpha^\nu\in\langle \eta\rangle$.
Let $A=\diag[\alpha\varepsilon,\alpha,1]$, $C=\diag[1,1,\eta]$
and $E=[e_3,e_1,e_2]$. Clearly
$(\diag[\eta,\eta,1])=(C)^{-1}$.\\
 3-1) The case where $G\cong K_1$. We may assume $K_1=G$. Suppose $\nu=\mu=9$. Since ${\bf{V}}(f)$ is non-singular and $f_{C^{-1}}\sim f$,
$f$ has the form $f_9(x,y)+cz^9$ with $cf_9 \not=0$. Let $f_9(x,y)=\sum_{i=0}^9a_ix^iy^{9-i}$. The
condition $f_{A^{-1}}\sim f$ yields $f_{A^{-1}}=f$, hence $(\alpha^9\varepsilon^{i}-1)a_i=0$
for all $i$. Let $\alpha^9=\eta^\ell=\varepsilon^j$ ($j,\ \ell\in [0,8]$). Unless $j=0$, ${\bf{V}}(f)$ is
singular at $(1,0,0)$ or $(0,1,0)$. So $\alpha^9=1$, and $f=ax^9+by^9+cz^9$ with $abc\not=0$.
Suppose $[\nu,\mu]\not=[9,9]$, namely $[\nu,\mu]\in
\{[1,81],[3,27],[27,3],[81,1]\}$. Then ${\bf V}(f)$ is singular by Proposition \ref{lem 1.14}. Thus
$\nu=\mu=9$, $\alpha^9=1$ and ${\bf{V}}(f)$ is projectively equivalent to the Fermat curve of degree 9. \\
3-2) The case where $G\cong K_2$. We may assume $K_2=G$. Observe that $\nu\mu=27$. If $[\nu,\mu]=[1,27]$, then
${\bf V}(f)$ is singular by Proposition \ref{lem 1.14} (1). Even if $[\nu,\mu]=[27,1]$, ${\bf V}(f)$ is singular,
for we may assume $A=A'=\diag[\varepsilon,\varepsilon^\ell,1]$
($\ell\in [1,26]$) by Lemma \ref{lem 1.5} (3) and we can show that $I(A')=\emptyset$.
Next let $[\nu,\mu]=[3, 9]$. Since ${{\bf V}}(f)$ is non-singular and $f_{C^{-1}}\sim f$,
 it follows that $f_{C^{-1}}=f$, hence $f=f_9(x,y)+cz^9$ with $cf_9\not=0$.
Now the condition $f_{A^{-1}}\sim f$ implies $f_{A^{-1}}=f$. Since $\alpha^3\in \langle \eta\rangle$ and
$\eta^3\in \{\varepsilon,\ \varepsilon^2\}$, we have $\alpha^9=\varepsilon^j$ ($j\in [0,2]$). Writing
$f_9=\sum_{i=0}^9 a_ix^iy^{9-i}$, we get $(\alpha^9\varepsilon^{i}-1)a_i=0$ for all $i\in
[0,9]$. Therefore, if $j=1$ or $j=2$, then ${\bf V}(f)$ is singular at $(1,0,0)$ or $(0,1,0)$, respectively.
If $j=0$, then $f=a_0y^9+a_3x^3y^6+a_6x^6y^3+a_9x^9+cz^9$. Since $f$ is $E$-invariant, $a_3=a_6=0$ so that
$f=ax^9+by^9+cz^9$ with $abc\not=0$.
Finally let $[\nu,\mu]=[9,3]$.
The non-singular ${{\bf V}}(f)$ with the property $f_{C^{-1}}\sim f$ must
satisfy $f_{C^{-1}}=f$, so that $f=f_9(x,y)+f_6(x,y)z^3+f_3(x,y)z^6+cz^9$ with $cf_9\not=0$.
Let $f_9=\sum_{i=0}^9a_ix^iy^{9-i}$, $f_6=\sum_{i=0}^6b_ix^iy^{6-i}$ and $f_3=\sum_{i=0}^3
c_ix^iy^{3-i}$. As $\alpha^9\in \langle \eta\rangle$ and
$\eta\in\{\varepsilon^3,\varepsilon^6\}$, we have $\alpha^9=\varepsilon^{3j}$ ($j\in [0,2]$).
Since $c\not=0$, hence
 $f_{A^{-1}}=f$, we obtain
$$\alpha^9f_9(\varepsilon x,y)=f_9(x,y),\,
\alpha^6f_6(\varepsilon x,y)=f_6(x,y),\, \alpha^3f_3(\varepsilon x,y)=f_3(x,y).$$
The condition $\alpha^9f_9(\varepsilon x,y)=f_9(x,y)$, that is,
$(\varepsilon^{3j+i}-1)a_i=0$ ($i\in [0,9]$), defines two polynomials
$f=a_6x^6y^3+f_6(x,y)z^3+f_3(x,y)z^6+cz^9$ and $f=a_3x^3y^6+f_6(x,y)z^3+f_3(x,y)z^6+cz^9$ which
give rise to the curve ${{\bf V}}(f)$ singular at $(1,0,0)$, provided $j\in \{1,2\}$. Let $j=0$. Then $\alpha^9=1$, hence
$\alpha=\varepsilon^\ell$ for some $\ell\in [0,8]$. Now, $f_9=ax^9+by^9$ by the condition
$\alpha^9f_9(\varepsilon x,y)=f_9(x,y)$. Obviously, the conditions
$\alpha^6f_6(\varepsilon x,y)=f_6(x,y)$ and $\alpha^3f_3(\varepsilon x,y)=f_3(x,y)$ are
equivalent to
\[
(\varepsilon^{6\ell+i}-1)b_i=0\ \ (i\in [0,6]),\hspace{10mm}(\varepsilon^{3\ell+i}-1)c_i=0\ \ (i\in [0,3]).
\]
Observe that $b_i=0$ unless $i\in \{0,3,6\}$ and that $c_i=0$ unless $i\in\{0,3\}$. Moreover, if
$\ell\equiv \ell'\ (\modd\ 3)$, then $\ell$ and $\ell'$ pose the same condition on $i$. Besides
$f_{E^{-1}}\sim f$. Suppose $\ell=0$. Since $f_6=b_0y^6$ and $f_3=c_0y^3$, the  $(E)$-invariant curve
 ${\bf{V}}(f)$ has equation $f=ax^9+by^9+cz^9=0$. Suppose $\ell=1$. Since $f_6=b_3x^3y^3$ and $f_3=0$,
the $(E)$-invariant curve ${\bf{V}}(f)$  has equation $f=ax^9+by^9+cz^9+ex^3y^3z^3=0$.  Suppose $\ell=2$.
Since $f_6=b_6x^6$ and $f_3=c_3x^3$,  the  $(E)$-invariant curve ${\bf{V}}(f)$  has equation
$f=ax^9+by^9+cz^9=0$. Thus, if $\ell\in \{0,2,3,5,6,8\}$, then $f=ax^9+by^9+cz^9$, and if
$\ell\in\{1,4,7\}$, then $f=ax^9+by^9+cz^9+e x^3y^3z^3$. Clearly both polynomials give rise to
singular ${\bf{V}}(f)$  if $abc=0$. If $abce\not=0$, ${\bf{V}}(f)$ is projectively equivalent to
${\bf{V}}(f')$, where $f'=x^9+y^9+z^9+\lambda x^3y^3z^3$ for some $\lambda\in k^*$. Such a curve
is non-singular if and only if $\lambda^3+27\not=0$ by Lemma\ref{lem 1.12}. \end{proof}

\begin{corollary} \label{cor 3.6}
Theorem \ref{th 3.1} holds for $d=9$.
\end{corollary}

\begin{proof} Note that $g'=g-1=d(d-3)/2=3^3$. The possible values of
$|\Aut({\bf{V}}(f))|\geq 6d^2=2\cdot 3^5$ are
$$84g', 48g', 40g', 36g', 30g', 24g', 21g',20g',56g'/3, 18g'=6d^2.$$
By Lemma \ref{lem 3.5} (1), there exists no non-singular curve ${\bf{V}}(f)$ of degree $9$ such that
$|\Aut({\bf{V}}(f))|\in \{40g',\ 30g',\ 20g'\}$ . The
remaining values are divisible by $3^4$ except for $56g'/3$, which is divisible by $8$. First
we will show that there is no non-singular curve ${\bf{V}}(f)$ of degree $d$ satisfying
$|\Aut({\bf{V}}(f))|=56g'/3$. Let $G$ be a subgroup of $PGL(3,k)$ such that $|G|=8$, and
assume that a plane curve ${\bf{V}}(f)$  of degree $d$ is
$G$-invariant. Clearly $8\not\in\{36,60,72,168,216,360\}$, and $8$ is not divisible by $3$.
Therefore we may assume that $G$ is the group described in Proposition \ref{pr 2.2} (2) or (3),
even if $G$ fixes a line.
Let $A=\diag[\alpha \varepsilon,\alpha,1]$,
$C=\diag[1,1,\eta]$ and $B=[\beta e_2,\beta e_1,e_3]$ with $\nu=\ord(\varepsilon)$ and
$\mu=\ord(\eta)$. We will discuss the former case $H\cong \D{2\nu}$
first. Using the notation in Proposition \ref{pr 2.2}, we have $8=|G'|=2\nu\mu$. Since $\nu\geq 2$,
$[\nu,\mu]$ is equal to $[2,2]$ or $[4,1]$. Assume $[\nu,\mu]=[2,2]$.
If ${\bf{V}}(f)$ is non-singular, we arrive at a contradiction as follows. Since
$f_{C^{-1}}\sim f$, we have $f_{C^{-1}}= f$, for otherwise $z$ divides $f$. Hence
$$f=f_9(x,y)+f_7(x,y)z^2 + f_5(x,y)z^4 +f_3(x,y)z^6+f_1(x,y)z^8$$
with $f_1\not=0$. By the
condition $f_{B^{-1}}\sim f$ we have $f_1=ax+by$ with $ab\not=0$, while the condition
$f_{A^{-1}}\sim f$ implies $\alpha f_1(\varepsilon x,y)\sim f_1(x,y)$, a contradiction, for
$\varepsilon=-1$.
Assume $[\nu,\mu]=[4,1]$, hence $B=[e_2,e_1,e_3]$ and $\alpha^4=1$. There exist four pairs of $[A,B]$, namely $[A_i,B]$  where
$A_i=\diag[\varepsilon^i,\varepsilon^{i-1},1]$ ($i\in [1,4]$), but
it suffices to consider the two pairs $[A_1,B]$ and $[A_2,B]$,
for $BA_1B=A_4^{-1}$ and $BA_2B=A_3^{-1}$.
An $A_1$-invariant $f$ is singular, unless $f_{A_1^{-1}}=f$. Since $f_{B^{-1}}\sim f$, $z$ divides $f$.
An $A_2$-invariant $f$ is singular, unless $f_{A_2^{-1}}=\varepsilon^i f$ ($i\in I(A_2)=\{0,1,2\}$).
Since $f_{B^{-1}}\sim f$, $z$ divides $f$. In fact, if $f_{A_{2}^{-1}}=f$ for example, then $x^3y^6$ and $x^7y^2$ are the possible
monomials which may appear in $f$ and which is not divisible by $z$. Since $f$ is $B$-invariant, the two monomials do not
appear in $f$, so that $z$ divides $f$.
We deal with the latter case $H\cong \Z_{\nu}$. Observe that $|G|=\nu\mu=8$, where
$\mu=\ord(\eta)$. In this case, $G$ is isomorphic to $\Z_8$ or $\Z_2\times \Z_4$. In fact, if
$\mu\in \{1,8\}$, then $G\cong \Z_8$. Let $\mu=2$ so that $\alpha^4=\pm 1$. If $\alpha^4=-1$,
then $G\cong \Z_8$. If $\alpha^4=1$, then $\langle (A)\rangle\cong \Z_4$, hence $G\cong
\Z_2\times \Z_4$. Let $\mu=4$, which yields $\alpha^2\in \langle \eta\rangle$. If $\alpha^2\in
\{\eta,\ \eta^3\}$, then $\ord(A)=8$, so that $G\cong \Z_8$. If $\alpha^2\in \{\eta^0,\
\eta^2\}$, then $\alpha\in \langle \eta\rangle$, so that $G\cong\Z_2\times \Z_4$. By Lemma
\ref{lem 3.5} (1), (2) we see that $|\Aut({\bf{V}}(f))|\leq 6d^2<56g'/3$, provided the $G$-invariant ${\bf{V}}(f)$ is non-singular.

Let $G$ be an automorphism group of a non-singular plane algebraic curve of degree $d=9$ with  $|G|\geq 6d^2$. Then we have
shown $|G|=6d^2$, which is a multiple of $3^4$. So Theorem \ref{th 3.1} holds for $d=9$ by  Lemma
\ref{lem 1.13}, Lemma \ref{lem 1.12}, and Lemma \ref{lem 3.5} (3).\end{proof}

\begin{lemma} \label{lem 3.7} Let ${\bf{V}}(f)$ be a plane algebraic curve of degree $10$.
\begin{itemize}
\item[\rm(1)] If ${\bf{V}}(f)$ is either $\Z_{7}$-invariant, or $\Z_{13}$-invariant or
$\Z_{25}$-invariant, then ${\bf{V}}(f)$ is singular.
\item[\rm(2)]If ${\bf{V}}(f)$ is non-singular and $\Z_{5}\times \Z_{5}$-invariant with $|\Aut({\bf{V}}(f))|\geq
6\cdot 10^2$, then ${\bf{V}}(f)$ is projectively equivalent to the Fermat curve of equation $x^{10}+y^{10}+z^{10}=0$.
\end{itemize}
\end{lemma}

\begin{proof} (1) Let $A=\diag[\varepsilon,\varepsilon^j,1]$ with
$\ord(\varepsilon)=\nu$ and $j\in [0,\nu-1]$, for $\nu\in \{7,13,25\}$.
We may assume $\Z_\nu=\langle (A)\rangle$ by Lemma \ref{lem 1.5} (3).  Then $I_x(A)=\{10,\
9+j,\ 9\}$, $I_y(A)= \{10j,\ 9j+1,\ 9j\}$, and $I_z(A)=\{0,\ 1,\ j\}$. These are subsets of
$\Z/\nu\Z$. It is easy to show that $I(A)=\emptyset$, hence any $\Z_\nu$-invariant plane curve
${\bf{V}}(f)$ of degree $10$ has a singular point due to
Lemma \ref{lem 1.6}. Let $\nu=25$, for example. Assume $0\in I(A)$. Since $0\in I_x$, we have $j=16$. Now,
$I_y(A)=\{10,19,20\}$ so that $0\not\in I_y(A)$, a contradiction. Similarly $1\not\in I(A)$.
Assume $j\in I(A)$. Since $j\in I_x(A)$, either $j=9$ or $j=10$. In any case we can verify $j\not\in I_y(A)$.
The remaining cases $\nu=7$ and $\nu=13$ can be dealt with similarly.

  (2) Let $B=\diag[1,1,\vare]$ and $C=\diag[1,\vare,1]$, where
$\ord(\vare)=5$. We may assume that $\Z_{5}\times \Z_{5}=\langle(B),(C) \rangle$ by Lemma
\ref{lem 1.7}. Note that $f_{B^{-1}}=\vare^i f$ and $f_{C^{-1}}=\vare^j f$ for
some $i,j\in [0,4]$. If either $i\in [1,4]$ or $j\in [1,4]$, then one can easily see that
either $z$ or $y$ divides $f$, respectively.  Assume $i=j=0$ so that
$$f=f_{10}(x,y)+f_5(x,y)z^5+f_0z^{10}=a_0x^{10}+a_1x^5y^5+a_2y^{10}+(b_0x^5+b_1y^5)z^5+c_0z^{10}.$$
Since ${\bf{V}}(f)$ is non-singular, we may assume
$f=x^{10}+x^5(F_1y^5+F_2z^5)+g(y,z)$, where $g=y^{10}+F_3y^5z^5+z^{10}$. It remains to show
that $F_1=F_2=F_3=0$. By computation we obtain $\Hess(f)=5^3x^3y^3z^3h(x,y,z)$, where
$$h=
x^{15}18\cdot 16F_1F_2+\cdots +x^0(36)(F_1y^5+F_2z^5)\{8F_3y^{10}+(36-F_3^2)y^5z^5
+8F_3z^{10}\}.$$
Let ${\cal L}_r$ be the set of linear factors of $\Hess(f)$ of multiplicity
$r$. Suppose $F_1F_2\not=0$. Then $(x)\in {\cal L}_3$, for $x$ does not divide $h$. Similarly $\{(y),\ (z)\}\subset {\cal L}_3$.
In particular $\Aut({\bf{V}}(f))$ acts on ${\cal L}_3$ so that
$$|\Aut({\bf{V}}(f))|/|\Aut({\bf{V}}(f))_x|\leq |{\cal L}_3|.$$
Let $(T)\in \Aut({\bf{V}}(f))_x$ with $T=[t_{ij}]$ such that $t_{11}=1$, and $\ t_{12}=t_{13}=0$. Let
$M=[t_{ij}]$ ($i,j\in [2,3]$), which belongs to $GL(2,k)$. Since $f_{T^{-1}}\sim f$, we get $t_{21}=t_{31}=0$,
$f_{T^{-1}}=f$, hence $(F_1y^5+F_2z^5)_{M^{-1}}=F_1y^5+F_2z^5$ and $g_{M^{-1}}=g$.
Since $\Hess(g)=225y^3z^3\{8F_3(y^{10}+z^{10})+(36-F_3^2)y^5z^5\}$, $M$
is either diagonal or skew-diagonal, for $g_{M^{-1}}=g$ implies $(yz)_{M^{-1}}\sim yz$. Thus if $T$ satisfies $(T)\in
\Aut({\bf{V}}(f))_x$, then either $t_{21}=t_{31}=t_{23}=t_{32}=0$, $t_{22}^{5}=t_{33}^{5}=1$ or
$t_{21}=t_{31}=t_{22}=t_{33}=0$, $t_{23}^{5}=F_2/F_1$, $t_{32}^{5}=F_1/F_2$, for $(F_1y^5+F_2z^5)_{M^{-1}}=F_1y^5+F_2z^5$.
Therefore $|\Aut({\bf{V}}(f))_x|\leq 2\cdot5^2$, so that $|\Aut({\bf{V}}(f))|\leq 50|{\cal L}_3|<600$, for $|{\cal
L}_3|\leq 3+\deg h/3=8$. Therefore $F_1F_2=0$. Similarly $F_1F_3=F_2F_3=0$. Hence only
one $F_i$ can be nonvanishing, provided $F_i\not=0$ for some $i$. We may assume $F_1= F_2=0$ and
$F_3\not=0$ without loss of generality. Then $\Hess(f)=90x^8\Hess(g)$. Observe that $\Hess(g)$ has no linear factors of
multiplicity 3 except for $y$ and $z$. In particular if $(T)\in \Aut({\bf{V}}(f))_x$ with $t_{11}=1$ and
$t_{12}=t_{13}=0$, then $M\in GL(2,k)$ defined as above from $T$ must be either diagonal or
skew-diagonal, because $g_{M^{-1}}\sim g$ implies $(yz)_{M^{-1}}\sim yz$. Besides
$(x)\in {\cal L}_8$, $|{\cal L}_8|=1$, and $\Aut({\bf{V}}(f))$ acts on ${\cal L}_8$. Now,
$|\Aut({\bf{V}}(f))_x|\leq
200$. In fact, $T$ satisfies $t_{11}=1,\ t_{12}=t_{13}=t_{21}=t_{31}=0$, and $M$ is either
diagonal or skew-diagonal. Since $f_{T^{-1}}\sim f$ implies $f_{T^{-1}}=f$,
 $t_{22}^{10}=t_{33}^{10}=1$ if $M$ is diagonal. If $M$ is skew-diagonal, then $t_{32}^{10}=t_{23}^{10}=1$. Thus
$|\Aut({\bf{V}}(f))_x|\leq 200$. Therefore $|\Aut({\bf{V}}(f))|\leq |\Aut({\bf{V}}(f))_x||{\cal L}_8|\leq 200<600$. Thus
all $F_i$ ($i=1,2,3$) must vanish.  \end{proof}


\begin{corollary} \label{cor 3.8}
Theorem \ref{th 3.1} holds for $d=10$.
\end{corollary}

\begin{proof}
Now, $g'=d(d-3)/2=5\cdot 7$. The possible values of $|\Aut({\bf{V}}(f))|> 600=(17+1/7)g'$ are $$84g', 48g',
36g', 30g', 12\cdot 11g'/5, 24g', 12\cdot 13g'/7, 21g', 20g', 96g'/5, 18g'.$$ They are multiples of $7$ or $13$.
We see that $|\Aut({\bf{V}}(f))|\leq 600=2^3\cdot 3 \cdot
5^2$ by Lemma \ref{lem 3.7} (1). A subgroup $G$ of $PGL(3,k)$ of order $5^2$ is isomorphic
to $\Z_{25}$ or $\Z_5\times \Z_5$, see \cite{hal}. So Lemma \ref{lem 3.7} and Lemma \ref{lem
1.13} imply Theorem \ref{th 3.1} for $d=10$. \end{proof}

\begin{remark} For curves with only ordinary flexes, Corollary \ref{cor 3.8} also follows from \cite[Section 3.6]{aluffi}.
\end{remark}

\noindent
{\bf Acknowledgement.}
The author would like to express her sincere thanks to Hitoshi Kaneta for his valuable comments and suggestions.

\end{document}